%% file: arxiv.tex
\documentclass[a4paper, 12pt]{article}

\usepackage[english]{babel}
\usepackage[utf8]{inputenc}
\usepackage[pdftex]{graphicx}
\usepackage[left=1in, right = 1in]{geometry}
\usepackage{array}
\usepackage{authblk}

\newcommand{\subD}{{\textsc{SubD}} }
\newcommand{\subB}{{\textsc{SubB}} }
\newcommand{\sepD}{{\textsc{SepD}} }
\newcommand{\sepB}{{\textsc{SepB}} }
\newcommand{\singlecut}{{\textsc{U-B\&BC}} }
\newcommand{\multicut}{{\textsc{A-B\&BC}} }
\newcommand{\localbranching}[2]{{\textsc{A-B\&BC+LB-{#1}-{#2}}} }
\newcommand{\branchandcut}{{\textsc{B\&C}} }
\newcommand{\greedy}{{\textsc{Greedy}} }
\newcommand{\branching}{{\textsc{SepB}} }
\newcommand{\disjunctive}{{\textsc{SepD}} }
\newcommand{\subdirect}{{\textsc{SubD}} }
\newcommand{\subbenders}{{\textsc{SubB}} }

\usepackage{appendix}
\usepackage{color,soul}
\usepackage{xcolor}

\usepackage{authblk}
\usepackage{amsthm,amsmath,amssymb}
\usepackage{algorithm}
\usepackage{algorithmic}
\usepackage{hyperref}
\usepackage{graphicx}
\usepackage{comment}
\usepackage{array}
\usepackage{rotating, graphicx}
\usepackage{tikz}
\usepackage{enumitem}

\renewcommand{\brace}[1]{\left\{ #1 \right\}}
\renewcommand{\epsilon}{\varepsilon}

\newcolumntype{L}[1]{>{\raggedright\let\newline\\\arraybackslash\hspace{0pt}}m{#1}}
\newcolumntype{C}[1]{>{\centering\let\newline\\\arraybackslash\hspace{0pt}}m{#1}}
\newcolumntype{R}[1]{>{\raggedleft\let\newline\\\arraybackslash\hspace{0pt}}m{#1}}

\usepackage[framemethod=TikZ]{mdframed}

\newtheorem{corollary}{Corollary}[section]

\newcounter{boxedproposition}[section]
\renewcommand{\theboxedproposition}{\arabic{section}.\arabic{boxedproposition}}
\newenvironment{boxedproposition}[0]{%
\refstepcounter{boxedproposition}%
\mdfsetup{%
frametitle={%
\tikz[baseline=(current bounding box.east),outer sep=0pt]
\node[anchor=east,rectangle,fill=blue!20]
{\strut Proposition~\theboxedproposition};}}%
\mdfsetup{innertopmargin=10pt,linecolor=blue!20,%
linewidth=2pt,topline=true,%
frametitleaboveskip=\dimexpr-\ht\strutbox\relax
}
\begin{mdframed}[]\relax%
}{\end{mdframed}}

\usepackage{algorithm}
\usepackage{algorithmic}
\usepackage{booktabs}

\usepackage{natbib}
 \bibpunct[, ]{(}{)}{,}{a}{}{,}%
 \def\newblock{\ }%
\begin{document}
\title{Accelerated Benders Decomposition and Local Branching for Dynamic Maximum Covering Location Problems}

\author[1]{Steven Lamontagne}
\author[1]{Margarida Carvalho}
\author[2]{Ribal Atallah}

\affil[1]{CIRRELT and D\'epartement d'informatique et de recherche op\'erationnelle, Universit\'e de Montr\'eal} 
\affil[2]{Institut de Recherche d'Hydro-Qu\'ebec}

\maketitle
\begin{abstract}
   The maximum covering location problem (MCLP) is a key problem in facility location, with many applications and variants. One such variant is the dynamic (or multi-period) MCLP, which considers the installation of facilities across multiple time periods. To the best of our knowledge, no exact solution method has been proposed to tackle large-scale instances of this problem. To that end, in this work, we expand upon the current state-of-the-art, branch-and-Benders-cut solution method in the static case, discussing several acceleration techniques. Additionally, we propose a specialised local branching scheme, which uses a novel distance metric in its definition of subproblems and features a new method for efficiently and exactly solving the subproblems. These methods are then compared with extensive computational experiments.
\end{abstract}

\section{Introduction}
A classic problem in operations research is the optimal location of facilities according to different aspects such as users' preferences and installation cost. Within this category of problems is the maximum covering location problem (MCLP), attributed to~\cite{Church1974}. In the MCLP, each facility covers the users within a certain radius and, due to limited resources, it is not possible to open every facility. Thus, a decision maker must select a subset of facilities to open, with the goal of maximising the total coverage. Due to its simplicity and versatility, the MCLP has been used in a wide range of applications, including emergency services location~\citep{Gendreau2001, Degel2015, Nelas2020}, healthcare services~\citep{Bagherinejad2018, Alizadeh2021}, safety camera positioning~\citep{Dellolmo2014, Han2019}, ecological monitoring or conservation~\citep{Church1996, MartinFores2021}, bike sharing~\citep{Muren2020}, and disaster relief~\citep{Zhang2017b, Iloglu2020, Yang2020}. For a review of the MCLP and its applications, we refer to \cite{Murray2016}.

To include more complex interactions and restrictions, several variants of the MCLP have been developed. The one we consider is the dynamic (or multi-period) MCLP, where the decision maker conducts facility planning over a long time horizon, divided into discrete time periods~\citep{Schilling1980, Gunawardane1982}. The intrinsic consideration of the time horizon is vital in many situations, such as when the demand to be covered varies across time~\citep{Porras2019}, when planning infrastructure that will persist throughout future planning periods~\citep{Gunawardane1982, Lamontagne2022}, or in real-time operations of emergency services where exact positioning is important~\citep{Gendreau2001}. Throughout the time horizon, open facilities may be forced to remain open for the duration of the time horizon~\citep{Lamontagne2022}, may be allowed to change location with or without cost~\citep{Marin2018}, or a subset may be required to relocated each time period~\citep{Dellolmo2014}. 

To the best of our knowledge, no exact method exists for solving large-scale dynamic MCLPs. More specifically, the only exact method employed in the literature for the dynamic MCLP is branch-and-bound, which has successfully solved instances with $300$ users, $300$ facilities, $9$ time periods~\citep{Dellolmo2014},  and $49,905$ users, $60$ facilities, $4$ time periods~\citep{Lamontagne2022}. Meanwhile, in terms of heuristic methods, \cite{Porras2019} used simulated annealing to solve instances with $547$ users, $547$ facilities, $4$ time periods, while \cite{Lamontagne2022} used a greedy method to solve instances with $649,605$ users, $180$ facilities, $4$ time periods. This compares with specialised methods for the static MCLP, which have successfully solved instances with $818$ users, $818$ facilities~\citep{Pereira2007}, and $15,000,000$ users, $100$ facilities~\citep{Cordeau2019}. Most notably, the method proposed in \cite{Cordeau2019} is specifically designed for problems in which there are significantly more users than facilities. This situation can occur when using a very fine discretisation of a continuous demand~\citep[e.g. ][]{Lamontagne2022} or when creating many different scenarios to model uncertainty in the problem~\citep{Daskin1983, Berman2013, Vatsa2016, Nelas2020}. These categories of use cases support the need for a specialised method for the dynamic MCLP which can handle large-scale instances.

In this work, we present several contributions to the literature on the dynamic MCLP. First, we extend the branch-and-Benders-cut approach for the static MCLP of \cite{Cordeau2019} to the dynamic case. Second, we detail a suite of acceleration techniques tailored for the dynamic MCLP, which can be selected based on the structure of the application. This includes an intuitive and effective multi-cut generation technique, an efficient Pareto-optimal cut generation technique with a closed-form solution, a partial Benders decomposition strategy which leverages the problem structure, and a Benders dual decomposition approach which can further strengthen cuts for fractional solutions.  Third, we present a specialized local branching method embedded within the branch-and-Benders cut framework. This method makes use of a new distance metric for defining subproblems, based on the problem structure. We then present a novel subproblem solution method, which can efficiently solve the local branching subproblems in an exact and proven-valid manner. Fourth, we present extensive computational results comparing our proposed methods. These experiments are carried out with instances based on real data from an electric vehicle charging station placement case study, underlying the practical interest of the dynamic MCLP. Our results validate the suitability of our methodological contributions to provide high-quality solutions to instances involving a large number of users, candidate locations and a planning horizon. In particular, our methods can produce solutions with a higher demand covered compared to the greedy heuristic in \cite{Lamontagne2022}, in addition to providing performance guarantees via optimality gaps.

The rest of this paper is organized as follows: Section~\ref{SectionLiteratureReview} discusses the literature about solution methods for the static and dynamic MCLP. Section~\ref{SectionProblemFormulation} presents the general model for the dynamic MCLP. Section~\ref{SectionBendersDecomposition} is dedicated to the branch-and-Benders-cut methods, extending the work in \cite{Cordeau2019} and implementing improvement methods to the general framework. Section~\ref{SectionLocalBranching} then presents the local branching method, based on the work in \cite{Rei2009}, with special consideration for the distance metric and the branching scheme. Finally, Section~\ref{SectionResults} provides our computational results, while Section~\ref{SectionConclusion} concludes our work.

\section{Literature Review}
\label{SectionLiteratureReview}

Despite its long history, few exact solution methods have been developped for the MCLP, while heuristic methods are most commonly used~\citep{Murray2016}. We provide a summary of solution methods for both the static and dynamic MCLP in Table~\ref{TableMaximumCoveringSolutionMethods}. 

The most common exact method uses standard branch-and-bound techniques, in conjunction with off-the-shelf mixed-integer linear programming solvers. This method is commonly used for small-scale instances, or as a benchmark approach for heuristics. However, it has been repeatedly noted that this is insufficient for solving large-scale instances~\citep[see, e.g., ][]{Revelle2008, Zarandi2013, Cordeau2019, Lamontagne2022}. In addition, three other exact methods have been presented for the static MCLP. In \cite{Downs1996}, a lower bound is generated via a greedy procedure, similar to the one presented in \cite{Church1974}, while an upper bound is generated via Lagrangian relaxation of the covering constraints. Both are combined within a branch-and-bound framework to ensure optimality. In \cite{Pereira2007}, the MCLP is reformulated as a p-median problem, with a column generation procedure used to solve the resulting model. In the procedure, stabilisation techniques are proposed to help improve the convergence rate, and limits are placed on both the number of columns generated and the total number of iterations. In \cite{Cordeau2019}, the variables associated with the coverage of users are projected out in the Benders subproblem. An analytic solution is then given for the solution of the Benders dual subproblem, allowing for rapid generation of Benders cuts.  This process is embedded within a branch-and-Benders-cut framework, allowing for large-scale problems to be solved efficiently.

In terms of heuristics, we note a wide variety of methods. However, two approaches which are used repeatedly include greedy methods and simulated annealing. In the greedy method~\citep[e.g. ][]{Church1974}, a solution is constructed by repeatedly opening a facility, each time selecting the one which causes the largest increase in coverage. In simulated annealing~\citep[e.g. ][]{Murray1996}, the process simulates the cooling of metals, with a temperature parameter that varies across the solving process and which controls the mutability of our current solution (e.g. swapping an open facility for a closed one). In \cite{Xia2009}, these two heuristics (along with greedy plus substitution, genetic algorithm and tabu search) were compared, with the simulated annealing method finding solutions of the best quality. By comparison, the greedy method solved instances faster, with solutions of value within a 1\% relative gap of those of simulated annealing. 

Our work contributes to the literature on exact methods for dynamic MCLPs, which, to our knowledge, is limited to branch-and-bound techniques. These have proved inadequate for tackling large-scale instances.

\begin{table}[]
    \centering
    \resizebox{0.9\textwidth}{!}{
    \begin{tabular}{l r r r}
    \hline
    Article & Variant & Methods proposed & Exact versus heuristic \\
    \hline
        \cite{Church1974} &  Static & Greedy, branch-and-bound & Both \\
        \cite{Schilling1980} & Dynamic & Weighing method & Heuristic   \\
        \cite{Gunawardane1982} & Dynamic & Branch-and-bound & Exact  \\
        \cite{Downs1996} &  Static & Greedy + Lagrangian relaxation + Branch-and-bound & Exact \\
        \cite{Galvao1996} &  Static & Greedy + Lagrangian relaxation & Heuristic \\
        \cite{Murray1996} & Static & Greedy, simulated annealing & Heuristic \\
        \cite{AdensoDiaz1997} & Static & Tabu search & Heuristic \\
        \cite{Resende1998} & Static & GRASP & Heuristic \\
        \cite{Galvao2000} & Static & Lagrangian relaxation, surrogate relaxation & Heuristic \\
        \cite{Arakaki2001} & Static & Genetic algorithm & Heuristic \\
        \cite{Gendreau2001} & Dynamic & Parallel tabu search & Heuristic \\
        \cite{Pereira2007} & Static & Column generation & Exact \\
        \cite{Revelle2008} & Static & Heuristic concentration & Heuristic\\
        \cite{Xia2009} & Static & Greedy, simulated annealing, genetic algorithm, tabu search & Heuristic \\
        \cite{Rodriguez2012} & Static & Iterative greedy & Heuristic \\
        \cite{Zarandi2013} & Dynamic & Simulated annealing & Heuristic \\
        \cite{Dellolmo2014} & Dynamic & Branch-and-bound & Exact \\
        \cite{Colombo2016} & Static & Variable neighbourhood search, heuristic concentration & Heuristic \\
        \cite{Calderin2017} & Dynamic & Simulated annealing, evolutionary algorithm & Heuristic \\
        \cite{Maximo2017} & Static & Intelligent-guided adaptive search & Heuristic \\
        \cite{Marin2018} & Dynamic & Lagrangian relaxation & Heuristic \\
        \cite{Cordeau2019} & Static & Branch-and-Benders-cut & Exact \\
        \cite{Porras2019} & Dynamic & Simulated annealing & Heuristic \\
        \cite{Lamontagne2022} & Dynamic & Greedy, GRASP, rolling horizon & Heuristic \\
        \hline
    \end{tabular}
    }
    \caption{Articles proposing solution methods for the static or dynamic MCLP. }
    \label{TableMaximumCoveringSolutionMethods}
\end{table}

\section{Problem Formulation}
\label{SectionProblemFormulation}
We present a general formulation for the dynamic MCLP. Let $T$ be the number of time periods, $J$ be the set of users and $I$ be the set of facilities. We use the index $1\leq t \leq T$ to denote a time period, the index $j \in J$ for a user, and the index $i \in I$ for a facility. We define the parameters $a_{ij}^t$ such that $a_{ij}^t = 1$ if facility $i$ covers user $j$ in period $t$, and $0$ otherwise. Let $d_j^t > 0$ denote the demand of user $j$ in time period $t$. We introduce a binary decision vector $x$ with entries $x^t_i$ such that  $x^t_i=1$ if and only if facility $i$ exists in period $t$. We also include in the model an auxiliary decision vector $z$ with entries $z^t_j$ such that $z^t_j=1$ if user $j$ is covered by some facility at period $t$, and $0$ otherwise. Then, the dynamic maximum covering location is formulated as 
\begin{subequations}
\begin{alignat}{3}
\operatorname{Maximise} & \sum_{t = 1}^T \sum_{j \in J} d_j^t z_j^t, 
\label{ModelMaximumCover:Objective} 
\\
\text{subject to }
& x \in \Omega,  
\label{ModelMaximumCover:Domain}
\\
& \sum_{i \in I}  a_{ij}^{t} x_i^t \geq z_j^t, && 1 \leq t \leq T, j \in J, 
\label{ModelMaximumCover:Covering}
\\
& x_{i}^t \in \brace{0, 1}, && 1 \leq t \leq T, i \in I, 
\label{ModelMaximumCover:BinaryX}
\\ 
& z_j^t \in \brace{0, 1} && 1 \leq t \leq T, j \in J.
\label{ModelMaximumCover:BinaryZ}
\end{alignat}
\label{ModelMaximumCover}
\end{subequations}
The objective function~\eqref{ModelMaximumCover:Objective} maximises the demand covered over all time periods. Constraint~\eqref{ModelMaximumCover:Domain} defines the feasible domain $\Omega$ for the variables $x$, which can be any polytope not involving the variables $z_j^t$. Typically, this includes either a cardinality constraint~\citep[e.g.][]{Calderin2017} or a knapsack constraint~\citep[e.g.][]{Cordeau2019} to restrict the number of facilities. The set $\Omega$ can also include, for instance, diversification constraints imposing that the set of facilities must change in each time period~\citep{Dellolmo2014}, or precedence constraints imposing that some facilities must be constructed before others~\citep{Lamontagne2022}. Constraints~\eqref{ModelMaximumCover:Covering} impose that a suitable facility must exist for our users to be considered covered. Constraints~\eqref{ModelMaximumCover:BinaryX} and~\eqref{ModelMaximumCover:BinaryZ} indicate that both sets of variables $x$ and $z$ must be binary. However, as noted in~\cite{Murray2016}, integrality on the variables $z_j^t$ can be relaxed in Constraints~\eqref{ModelMaximumCover:BinaryZ}.

For a given solution $x \in \Omega$, it is easy to calculate the resulting coverage or, in other words, the total demand covered by the set of facilities defined via $x$. Let $f(x)$ denote the coverage of solution $x \in \Omega$. Then, as described in~\cite{Revelle2008, Cordeau2019}, this value can be calculated  
as
\begin{equation}
\label{SolutionCoverage}
    f(x) = \sum_{t=1}^T \sum_{j \in J} \min \brace{1, \sum_{i \in I} a_{ij}^t x_i^t} d_j^t.
\end{equation}

\section{Benders Decomposition}
\label{SectionBendersDecomposition}

In this section, we present an accelerated branch-and-Benders-cut approach tailored for the dynamic MCLP. This begins by generalising the current state-of-the-art approach in the static MCLP, the branch-and-Benders-cut method proposed by \cite{Cordeau2019}, to the dynamic case. We then discuss acceleration techniques which aim to strengthen the formulation and accelerate convergence, applied within our dynamic context.

\subsection{Single Cut Benders Decomposition}
\label{SectionSingleCut}

We detail the development of the branch-and-Benders-cut method proposed in \cite{Cordeau2019}. While the process is nearly identical in the dynamic case as the static case, the development process is a prerequisite for the acceleration techniques we discuss below. As in \cite{Cordeau2019}, we define $J_s= \left\{ j \in J : \sum_{t=1}^T \sum_{i \in I} a_{ij}^t = 1 \right\}$ as the set of users which are covered by only one facility. 

Like in the classical Benders decomposition~\citep{Benders1962}, we project out the (continuous) $z$ variables of model~\eqref{ModelMaximumCover}. The value of the $z$ variables in the objective function is replaced with an auxiliary variable, $\theta$.  Then, in iteration $v$,  the main problem~(MP) can be written as
\begin{subequations}
	\begin{alignat}{5}
\operatorname{Maximise } & \, \, \theta, 
\\
\text{subject to }
& \eqref{ModelMaximumCover:Domain}, \eqref{ModelMaximumCover:BinaryX} \notag
\\ 
& \operatorname{Opt}_{\pi^r, \sigma^r} (x) \geq \theta, \quad && 1 \leq r \leq v, 
\label{SingleCut:BendersMain:OptimalityCuts}
\\
& \theta \geq 0,
\label{SingleCut:BendersMain:DomainTheta}
\end{alignat}
\label{SingleCut:BendersMain}
\end{subequations}
where $\pi^r$ and $\sigma^r$ are the optimal dual vectors associated with the Benders subproblem in iteration $r$. Since the Benders subproblem is feasible $\forall x \in \Omega$, it is not necessary to include feasibility cuts. 

We denote —now and throughout this paper— $\tilde{x}$ as the candidate solution in iteration $v$. Unless otherwise specified, it is assumed that this solution is integer feasible, i.e. that $\tilde{x} \in \Omega$ and $\tilde{x}_i^t \in \{0,1\}, 1 \leq t \leq T, i \in I$.  For each $1 \leq t \leq T$ and $j \in J$, let $I_j^t(x) = \sum_{i \in I} a_{ij}^t x_{i}^t$ denote the coverage of user $j$ in period $t$ by solution $x$. When there is no risk of confusion, we use $\tilde{I}_j^t = I_j^t(\tilde{x})$.

The Benders Primal Subproblem~(BPS) is then given by
\begin{subequations}
	\begin{alignat}{5}
\operatorname{Maximise} & \, \, \sum_{t = 1}^T \sum_{j \in J} d_j^t z_j^t,
\\ \text{subject to }
& z_j^t \leq \tilde{I}_j^t, && 1 \leq t \leq T, j \in J, 
\label{SingleCut:BendersSubproblemPrimal:Covering}
\\
& z_j^t \leq 1, && 1 \leq t \leq T, j \in J,
\label{SingleCut:BendersSubproblemPrimal:UpperBoundZ}
\\
& z_j^t \geq 0, && 1 \leq t \leq T, j \in J.
\label{SingleCut:BendersSubproblemPrimal:LowerBoundZ}
\end{alignat}
\label{SingleCut:BendersSubproblemPrimal}
\end{subequations}

Let $\pi_{j}^{t}$ and $\sigma_j^{t}$ denote the dual variables associated with Constraints~\eqref{SingleCut:BendersSubproblemPrimal:Covering} and~\eqref{SingleCut:BendersSubproblemPrimal:UpperBoundZ}. The Benders Dual Subproblem~(BDS) is then
\begin{subequations}
\begin{alignat}{5}
\operatorname{Minimise} & \sum_{t=1}^T \sum_{j \in J} \left(\tilde{I}_j^{t} \pi_{j}^{t}  + \sigma_j^{t}\right),  
\\ \text{subject to }
& \pi_{j}^{t} + \sigma_{j}^{t} \geq d_j^t, \, && 1 \leq t \leq T, j \in J, 
\\
& \pi_{j}^{t}, \sigma_{j}^{t} \geq 0, \, && 1 \leq t \leq T, j \in J.
\end{alignat}
\label{SingleCut:BendersSubproblemDual}
\end{subequations}
Similar to \cite{Cordeau2019}, this subproblem can be easily solved by inspection: If $\tilde{I}_j^t < 1$, then the optimal solutions are $\tilde{\pi}_j^t = d_j^t, \tilde{\sigma}_j^t = 0$. If $\tilde{I}_j^t > 1$, then $\tilde{\pi}_j^t = 0, \tilde{\sigma}_j^t = d_j^t$. If $\tilde{I}_j^t = 1$, then any solution $\tilde{\pi}, \tilde{\sigma} \geq 0$ such that $\tilde{\pi}_j^t + \tilde{\sigma}_j^t = d_j^t$ will be optimal.  The optimality cut associated with these values is then
$$
\sum_{t=1}^T \sum_{j \in J} \left( \left(\sum_{i \in I} a_{ij}^t x_i^t \right) \tilde{\pi}_{j}^{t}  + \tilde{\sigma}_j^{t} \right) \geq \theta.
$$
The resulting cut, reformulated as a function of the variables $x$, is presented in Proposition~\ref{PropositionCuts}. We set $\pi^{v+1}, \sigma^{v+1} = \tilde{\pi}, \tilde{\sigma}$, and the left-hand side of the cut forms the term $\operatorname{Opt}_{\pi^{v+1}, \sigma^{v+1}}$ in \eqref{SingleCut:BendersMain:OptimalityCuts}.  

\begin{boxedproposition}
\label{PropositionCuts}
The optimality cuts associated with the Benders subproblem take the form, 
\begin{equation}
\label{OptimalityCutSingle}
    \sum_{t = 1}^T \sum_{i \in I} \left( \sum_{j \in \Gamma^t(\tilde{x}^t)} d_j^t a_{ij}^t \right) x_i^t + \sum_{t=1}^T \sum_{j \in J\setminus \Gamma^t(\tilde{x}^t)} d_j^t \geq \theta 
\end{equation}
where the set $\Gamma^t(\tilde{x}^t)$ can be defined as any of the following expressions:
\begin{align}
\tag{B0} \label{SingleB0} 
\Gamma^t(\tilde{x}^t) & = \left\{ j \in J: \tilde{I}_j^t < 1 \right\}, \\
\tag{B1} \label{SingleB1} 
\Gamma^t(\tilde{x}^t) & = \left\{ j \in J\setminus J_s: \tilde{I}_j^t < 1 \right\} \cup \left\{ j \in J_s: \tilde{I}_j^t \leq  1 \right\}, \\
\tag{B2} \label{SingleB2} 
\Gamma^t(\tilde{x}^t) & = \left\{ j \in J: \tilde{I}_j^t \leq 1 \right\}.
\end{align}
\end{boxedproposition}
We remark that the names of the cuts~\eqref{SingleB0},~\eqref{SingleB1}, and~\eqref{SingleB2} correspond to those of the equivalent cuts in \cite{Cordeau2019}.

\begin{figure}
    \centering
\resizebox{0.5\textwidth}{!}{
\input{Image_SubproblemMethod}
}
    \caption{A simple example, with $T=1$, three facilities $\{i_1, i_2, i_3 \}$ and user coverage given by the numbers within the regions. Numbers within the intersection of facilities indicate overlapping coverage. }
    \label{ExampleOverlap}
\end{figure}
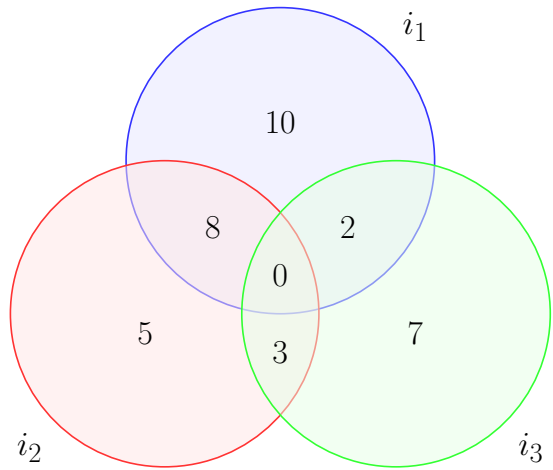

We note the following observations about these optimality cuts:
\begin{itemize}
    \item As in \cite{Cordeau2019}, the cuts~\eqref{SingleB0} are dominated by the  cuts~\eqref{SingleB1}. 
    \item By definition of the set $J_s$, the condition $\tilde{I}_j^t \leq  1$ will always be satisfied for $j \in J_s$.
    \item The optimality cuts~~\eqref{OptimalityCutSingle} are valid for both integer and fractional candidate solutions $\tilde{x}$. As described in \cite{Cordeau2019}, optimality cuts can be generated for fractional solutions to improve the convergence rate.  
    \item  When evaluated at solution $\tilde{x}$, the value of the left-hand side of the optimality cut~\eqref{OptimalityCutSingle} (and, thus, the corresponding objective value) matches the value of $\tilde{x}$ in Equation~\eqref{SolutionCoverage}. This derives directly from the definition of $\Gamma^t(\tilde{x}^t)$.  However, this evaluation is not correct (i.e. equal to $f(\tilde{x}$ in Equation~\eqref{SolutionCoverage}) in general if evaluated at a different solution. 
    
    We present a simple example of this in Figure~\ref{ExampleOverlap}. We drop the superscript $t$ as there is only one time period. If we generate the optimality cuts for the candidate solution $\tilde{x} = (1,0,0)$, i.e. facility $i_1$ is opened whereas facilities $i_2$ and $i_3$ are closed, we have the following optimality cuts
\begin{align*}
    \text{(B0): } \, & 0 x_{i_1} + 8 x_{i_2} + 10 x_{i_3} + 20 \geq \theta, \\
    \text{(B1): } \, & 10 x_{i_1} + 8 x_{i_2} + 10 x_{i_3} + 10 \geq \theta, \\
    \text{(B2): } \, & 20 x_{i_1} + 16 x_{i_2} + 12 x_{i_3} + 0 \geq \theta.
\end{align*}
All three types of cuts correctly identify the objective value of solution $\tilde{x}$ as $20$, matching the value in Equation~\eqref{SolutionCoverage}.  However, if we then use these optimality cuts to calculate the objective value of the solution $\bar{x} = (0, 0.75, 0.75)$, we obtain the values $33.5, 23.5,$ and $21$ (for B0, B1, and B2 cuts respectively), while Equation~\eqref{SolutionCoverage} gives $19.5$.
\end{itemize}

So far, we have presented a simple generalisation of the Benders procedure of \cite{Cordeau2019} to the dynamic MCLP. However, as we will see in Section~\ref{SectionResultsBaseline}, this method can exhibit slow convergence. We next explore improvements with respects to its convergence rate.

\subsection{Improvements to the Single-Cut Method}
\label{SectionImprovements}

In this section, we demonstrate and discuss the application of known techniques to our particular dynamic MCLP. These approaches aim to strengthen the Benders optimality cuts and improve the dual bound for the decomposition.   For the sake of simplicity, the acceleration techniques are presented separately in relation to the decomposition presented in Section~\ref{SectionSingleCut}. However, it is possible to incorporate multiple techniques simultaneously. For a discussion of acceleration techniques in the general context of Benders decomposition, we refer to \cite{Rahmaniani2017}.

There exist classes of improvement methods which heavily depend on the structure of $\Omega$, but which do not affect the structure of the Benders subproblems. Since they involve $\Omega$, their use is application-specific. However, since they do not impact the Benders subproblems, they can be applied to the methods in this work without affecting their validity. An example of such a technique would be a primal heuristic providing warmstart solutions. This can be used to obtain a good lower bound on the optimal objective value,  which otherwise can be difficult during the early iterations of the Benders decomposition method~\citep{Rahmaniani2017}. Another example would be valid cuts, which can be added to the main problem to tighten its linear programming (LP) relaxation. Examples for both techniques can be found in \cite{Santoso2005, Codato2006, Costa2012}.

The remainder of the techniques discussed in this section modify the Benders subproblems in some way, whether by changing the set of users $J$ (which, in turn, changes the amount of computational work required for the generation of cuts), modifying the frequency or type of cuts, or by directly providing reformulations for the BPS~\eqref{SingleCut:BendersSubproblemPrimal} or the BDS~\eqref{SingleCut:BendersSubproblemDual}. When an aspect of the technique depends on the structure of $\Omega$ (and thus may be application-specific), we include a discussion on this dependence. 

\paragraph{Preprocessing} 
Preprocessing techniques which aggregate users reduce the size of $J$ and, consequently, the number of $\tilde{I}_j^t$ to be calculated in each iteration of the Benders decomposition. A method for doing this such that the resulting problem is equivalent is given in \cite{Legault2022}. Their method aggregates users based on their coverage from facilities, i.e., if users $j$ and $j'$ are covered by exactly the same set of facilities then they can be aggregated. Heuristic methods also exist such as the ones by \cite{Dupacova2003} and \cite{Crainic2014}. However, we note that, due to the analytic solution for the Benders dual subproblem, the marginal effects of each user is negligible in terms of the cut generation time. 

\paragraph{Cut generation at fractional solutions} In Section~\ref{SectionSingleCut}, we observed that the Benders optimality cuts could be generated for both integer and fractional candidate solutions $\tilde{x}$. In \cite{Cordeau2019}, these cuts are generated at every node in the search tree, regardless if the corresponding solution is fractional or integer. However, in \cite{Botton2013}, it was demonstrated that this can lead to an increase to the solving time, due to the computational effort required to generate the cuts with marginal benefits. As such, it was proposed to generate cuts at fractional solutions only at the root node. 

\paragraph{Multi-cut method}
\label{SectionMulticut} 
Since the variables $z^t_{ij}$  are only involved in the covering Constraints~\eqref{ModelMaximumCover:Covering}, the Benders primal subproblem~\eqref{SingleCut:BendersSubproblemPrimal} can be separated per time period. This allows for more optimality cuts to be generated each iteration of the Benders decomposition, which accelerates convergence, as noted by \cite{Birge1988} and \cite{Contreras2011}. More specifically, we replace the auxiliary variable $\theta$ by the sum of new variables $\theta^t$ for $1\leq t \leq T$. Then, for each iteration $v$, we have the following MP:

\begin{subequations}
	\begin{alignat}{5}
\operatorname{Maximise } & \, \, \sum_{t=1}^t \theta^t, 
\\
\text{subject to }
& \eqref{ModelMaximumCover:Domain}, \eqref{ModelMaximumCover:BinaryX}
\\ 
& \operatorname{Opt}_{\pi^{rt}, \sigma^{rt}} (x^t) \geq \theta^t, \quad && 1 \leq t \leq T, 1 \leq r \leq v,
\label{MultiCut:BendersMain:OptimalityCuts}
\\
& \theta^t \geq 0,
\label{MultiCut:BendersMain:DomainTheta}
\end{alignat}
\end{subequations}
where, as before, $\pi^{rt}$ and $\sigma^{rt}$ are the optimal dual vectors associated with the Benders subproblem in time period $t$ and iteration $r$. The BPS~\eqref{SingleCut:BendersSubproblemPrimal} and the  BDS~\eqref{SingleCut:BendersSubproblemDual} as well as the process for deriving the Benders cuts is identical to the single-cut case, with the summations (and, as a consequence, indices) shifted from the subproblem to the main problem. The optimality cuts themselves are presented in Proposition~\ref{PropositionMulticut}.

\begin{boxedproposition}
\label{PropositionMulticut}
For each $1\leq t\leq T$, the optimality cuts associated with the Benders subproblem take the form, 
\begin{equation}
\label{OptimalityCutMulti}
    \sum_{i \in I} \left( \sum_{j \in \Gamma^t(\tilde{x}^t)} d_j^t a_{ij}^t \right) x_j^t + \sum_{j \in J\setminus \Gamma^t(\tilde{x}^t)} d_j^t \geq \theta^t .
\end{equation}
The sets $\Gamma^t(\tilde{x}^t)$ are as defined in Proposition~\ref{PropositionCuts}.
\end{boxedproposition}

While the Benders primal subproblem~\eqref{SingleCut:BendersSubproblemPrimal} can be separated per time period, it can also be separated by user. In fact, the preceding development can be done using auxiliary variables $\theta^t_j$ for $1\leq t \leq T, j \in J$ and Benders primal subproblems defined for each time period and user. However, as this method is designed for large-scale instances with many users, such a decomposition may be computationally infeasible. Notably, while the multi-cut generation process given in Proposition~\ref{PropositionMulticut} adds $|T|$ cuts in each iteration of the Benders decomposition, the equivalent version when separated by time period and by user would generated $|T| \times |J|$ cuts in each iteration. As such, the number of constraints in the model would grow extremely rapidly each iteration of the Benders decomposition, quickly becoming intractable.

\paragraph{Pareto-optimal cuts}
\label{SectionParetoOptimal}

The optimal solution to the Benders dual subproblem may not be unique, and the strength of the resulting Benders cut may vary depending on the solution selected. This phenomenon is the principle behind the B0, B1, and B2-type cuts introduced in Proposition~\ref{PropositionCuts}, with the strength of the resulting cuts discussed in depth in \cite{Cordeau2019}. In the context of general Benders decomposition, the selection of an optimal solution leading to a Benders cut is introduced in \cite{Magnanti1981}. The authors propose a model for finding the strongest possible cut based on a \emph{core point}, a point in the relative interior of the convex hull of $\Omega$. This technique was employed in \cite{Santoso2005, Contreras2011}, and discussed further in \cite{Papadakos2008}.

Let $\tilde{x}$ be the candidate solution from the main problem~\eqref{SingleCut:BendersMain} and  $\tilde{\pi}, \tilde{\sigma}$ be an associated optimal solution of the dual subproblem~\eqref{SingleCut:BendersSubproblemDual}. Let $\overset{c}{x}$ denote a core point of $\Omega$. \cite{Magnanti1981} define a Pareto-optimal cut with respect to a core point $\overset{c}{x}$ and a candidate solution $\tilde{x}$ as the one obtained by solving the following subproblem:
\begin{subequations}
\begin{alignat}{5}
(\overset{c}{\pi}, \overset{c}{\sigma}) \in \operatorname{argmin} & \sum_{t=1}^T \sum_{j \in J} \left( I(\overset{c}{x})_{j}^{t} \pi_{j}^{t}  + \sigma_j^{t} \right),  
\\ \text{subject to }
& \pi_{j}^{t} + \sigma_{j}^{t} \geq d_j^t, \, && \hspace{-0.5cm}1 \leq t \leq T, j \in J, 
\label{ParetoOptimal:MWSubproblem:DualFeas}
\\
& \sum_{t=1}^T \sum_{j \in J} I(\tilde{x})_j^{t} \pi_j^{t} + \sigma_j^{t} = \sum_{t=1}^T \sum_{j \in J}  I(\tilde{x})_j^{t} \tilde{\pi}_j^{t} + \tilde{\sigma}_j^{t},
\label{ParetoOptimal:MWSubproblem:Main}
\\
& \pi_{j}^{t}, \sigma_{j}^{t} \geq 0, \, && \hspace{-0.5cm}1 \leq t \leq T, j \in J,
\label{ParetoOptimal:MWSubproblem:NonNeg}
\end{alignat}
\label{ParetoOptimal:MWSubproblem}
\end{subequations}
where $(\overset{c}{\pi}, \overset{c}{\sigma})$  is called a Pareto-optimal point.

Due to the presence of Constraint~\eqref{ParetoOptimal:MWSubproblem:Main}, the Magnonti-Wong subproblem~\eqref{ParetoOptimal:MWSubproblem} is not quite identical to the Benders dual subproblem~\eqref{SingleCut:BendersSubproblemDual}. Nevertheless, we are still able to describe the determination of an optimal solution to subproblem~\eqref{ParetoOptimal:MWSubproblem}:

\begin{boxedproposition}
\label{TheoremParetoOptimal}
For a given candidate solution $\tilde{x}$ and core point $\overset{c}{x}$, the Pareto-optimal solutions $\overset{c}{\pi},\overset{c}{\sigma}$ for problem~\eqref{ParetoOptimal:MWSubproblem} can be calculated pointwise as follows: For $1 \leq t \leq T, j \in J$,
\begin{itemize}
    \item If $I(\tilde{x})_j^{t} < 1$, then $\overset{c}{\pi}_{j}^{t} = d_j^t, \overset{c}{\sigma}_{j}^{t} = 0$.
    \item If $I(\tilde{x})_j^{t} > 1$, then $\overset{c}{\pi}_{j}^{t} = 0, \overset{c}{\sigma}_{j}^{t} = d_j^t$.
    \item If $I(\tilde{x})_j^{t} = 1$ and $I(\overset{c}{x})_{j}^{t} < 1$, then $\overset{c}{\pi}_{j}^{t} = d_j^t, \overset{c}{\sigma}_{j}^{t} = 0$.
    \item If $I(\tilde{x})_j^{t} = 1$ and $I(\overset{c}{x})_{j}^{t} > 1$, then $\overset{c}{\pi}_{j}^{t} = 0, \overset{c}{\sigma}_{j}^{t} = d_j^t$.
    \item If $I(\tilde{x})_j^{t} = 1$ and $I(\overset{c}{x})_{j}^{t} = 1$, then  any value such that $\overset{c}{\pi}_{j}^{t} + \overset{c}{\sigma}_{j}^{t} = d_j^t$ is Pareto-optimal.   
\end{itemize}

\end{boxedproposition}
\begin{proof}
We first emphasise that the solution to the Magnanti-Wong subproblem~\eqref{ParetoOptimal:MWSubproblem} must also be an optimal solution to the BDS~\eqref{SingleCut:BendersSubproblemDual} for $\tilde{x}$, as specified in \cite{Magnanti1981}. This is enforced through Constraints~\eqref{ParetoOptimal:MWSubproblem:DualFeas} and~\eqref{ParetoOptimal:MWSubproblem:NonNeg} (which enforce feasibility), combined with Constraint~\eqref{ParetoOptimal:MWSubproblem:Main} (which enforces optimality). As described in Section~\ref{SectionSingleCut}, there is a set of optimal solutions to BDS~\eqref{SingleCut:BendersSubproblemDual} which can be analytically determined for each period $t$ and customer $j$ based on $I(\tilde{x})^t_j$. Hence, let us denote by $\Phi^t_j$ the set of optimal solution pairs $(\tilde{\pi}_{j}^{t}, \tilde{\sigma}_{j}^{t})$. Concretely, recall that if $I(\tilde{x})_j^{t} < 1$, then $\Phi_j^{t} = \brace{ (d_j^t, 0) }$. Likewise, $I(\tilde{x})_j^{t} > 1$ implies that $\Phi_j^{t} = \brace{ (0, d_j^t) }$. Finally, if $I(\tilde{x})_j^{t} = 1$ then $\Phi_j^{t} = \left\{ ({\pi}_{j}^{t}, {\sigma}_{j}^{t}) :  {\pi}_{j}^{t} + {\sigma}_{j}^{t} = d_j^t, {\pi}_{j}^{t} \geq 0,  {\sigma}_{j}^{t} \geq 0 \right\}$. 

In this way, we can equivalently formulate subproblem~\eqref{ParetoOptimal:MWSubproblem} with Constraint~\eqref{ParetoOptimal:MWSubproblem:Main} replaced by the following constraint:
$$
(\pi_j^t, \sigma_j^t) \in \Phi_j^t.
$$
This eliminates the linking constraint which involves all variables, reformulating the Magnonti-Wong subproblem~\eqref{ParetoOptimal:MWSubproblem} into the following:
\begin{subequations}
\begin{alignat}{5}
\operatorname{Minimise } \, & I(\overset{c}{x})_{j}^{t} \pi_{j}^{t}  + \sigma_j^{t},  
\\ \text{subject to }
& \eqref{ParetoOptimal:MWSubproblem:DualFeas},\eqref{ParetoOptimal:MWSubproblem:NonNeg} \notag
\\
& \pi_{j}^{t}, \sigma_{j}^{t} \in \Phi_j^{t}.
\label{MWSubproblemSeparated:NonNeg}
\end{alignat}
\label{MWSubproblemSeparated}
\end{subequations}

Now, it is easy to verify the feasibility and optimality of the solution described in the proposition statement. 
\end{proof}

\begin{corollary}
If the set $\brace{(t,j): 1 \leq t \leq T, j \in J, I(\tilde{x})_j^{t} =1}$ is empty, then the Benders optimality cut is unique.
\end{corollary}

In practice, Proposition~\ref{TheoremParetoOptimal} states that, when determining if $j \in \Gamma^t(\tilde{x}^t)$, we only need to consider the core point when $I_j^t(\tilde{x}) =1$. In those cases, rather than solely using information from $j$ and $\tilde{x}$ (as is the case for the cuts in Proposition~\ref{PropositionCuts}), we can instead look to the coverage in the core point to determine if $j \in \Gamma^t(\tilde{x}^t)$. If both $I_j^t(\tilde{x}) =1$ and $I(\overset{c}{x})_j^t =1$, then a different method is required. For example, it is possible to create analogous cuts to~\eqref{SingleB0},~\eqref{SingleB1}, or~\eqref{SingleB2}.

The calculation of Pareto-optimal points relies on finding a core point $\overset{c}{x}$. Depending on the domain $\Omega$, it may not be possible to find a core point analytically. However, one may employ the iterative method proposed in~\cite{Papadakos2008}, which starts with a feasible solution and takes the average with the candidate solution at every iteration to take the role of core point.

\paragraph{Benders Dual Decomposition}
\label{SectionBendersDualDecomposition}

In Benders Dual Decomposition, information from the main problem is added to the Benders subproblems to generate better quality cuts. This is done via allowing some (or potentially all) of the upper level variables to change based on the dual information. Following the work in \cite{Rahmaniani2020}, auxiliary variables $y_i^t$ are created by duplicating the decision variables $x_{i}^t$. These auxiliary variables are then forced to match a fractional or integer candidate solution $\tilde{x}_i^t$ in the Benders subproblem, resulting in the following subproblem:

\begin{subequations}
\begin{alignat}{3}
\operatorname{PSP}\left( \tilde{x} \right) = \operatorname{Maximise} & \, \, \sum_{t = 1}^T \sum_{j \in J} d_j^t z_j^t,
\\ 
\text{subject to } \,
& y_i^t = \tilde{x}_i^t, && 1 \leq t \leq T, i \in I,
\label{BDD:BendersSubproblem:MatchCandidate} 
\\
& z_j^t \leq \sum_{i \in I} a_{ij}^t y_i^t, && 1 \leq t \leq T, j \in J, 
\label{BDD:BendersSubproblem:Covering} 
\\
& y \in \Omega, 
\label{BDD:BendersSubproblem:DomainY} 
\\
& z_j^t \in [0,1], && 1 \leq t \leq T, j \in J.
\label{BDD:BendersSubproblem:DomainZ} 
\end{alignat}
\label{BDD:BendersSubproblem}
\end{subequations}

In its current form, the subproblem~\eqref{BDD:BendersSubproblem} is entirely equivalent to the original subproblem~\eqref{SingleCut:BendersSubproblemPrimal}. However, by applying Lagrangian relaxation on the Constraints~\eqref{BDD:BendersSubproblem:MatchCandidate}, for a given Lagrangian multiplier $\lambda_i^t$ for each $1\leq t \leq T$ and $i \in I$, we obtain the following Lagrangian subproblem:
\begin{subequations}
\begin{alignat}{3}
\operatorname{LSP1}\left( \tilde{x}, \lambda \right) = \operatorname{Maximise } \,  & \sum_{t=1}^T \sum_{i \in I} \left(d_j^t z_j^t - \lambda_{i}^t \left( y_{i}^t - \tilde{x}_{i}^t\right) \right), 
\label{BDD:BendersSubproblemLagrangianPartial:Objective} 
\\
\text{subject to }
& \eqref{BDD:BendersSubproblem:Covering}-\eqref{BDD:BendersSubproblem:DomainZ}.
\end{alignat}
\label{BDD:BendersSubproblemLagrangianPartial}
\end{subequations}
If we take  $\lambda_{i}^t = \tilde{\lambda}_{i}^t$ as the optimal dual variables associated with Constraints~\eqref{BDD:BendersSubproblem:MatchCandidate} after solving the subproblem~\eqref{BDD:BendersSubproblem} and $\bar{y}, \bar{z}$ as the optimal solutions of subproblem $\operatorname{LSP1}\left( \tilde{x}, \tilde{\lambda} \right)$, this results in the optimality cuts 
\begin{equation}
\label{BDD:StrengthenedCuts}
    \theta \leq \sum_{t=1}^t \sum_{j \in J} d_j^t \bar{z}_j^t - \sum_{t =1}^T \sum_{i \in I} \tilde{\lambda}_{i}^t \left( \bar{y}_{i}^t - x_{i}^t  \right).
\end{equation}
In \cite{Rahmaniani2020}, it was shown that the cuts~\eqref{BDD:StrengthenedCuts} are equivalent to the standard optimality cuts~\eqref{OptimalityCutSingle} when generated from an integer candidate solution $\tilde{x}$. On the other hand, the authors also show that the cuts~\eqref{BDD:StrengthenedCuts} are stronger than cuts~\eqref{OptimalityCutSingle} when generated from a fractional candidate solution.

If we further solve the full Lagrangian subproblem
\begin{equation}
\label{BDD:BendersSubproblemLagrangianFull}
  \operatorname{LSP2}\left( \tilde{x} \right) =   \operatorname{Minimise}_{\lambda} \left( \operatorname{LSP1}\left( \tilde{x}, \lambda \right) \right),
\end{equation}
with optimal solutions $\bar{\lambda}, \bar{y}, \bar{z}$, we obtain the following optimality cuts: 
\begin{equation}
\label{BDD:ExtraStrengthenedCuts}
    \theta \leq \sum_{t=1}^t \sum_{j \in J} d_j^t \bar{z}_j^t - \sum_{t =1}^T \sum_{i \in I} \bar{\lambda}_{i}^t \left( \bar{y}_{i}^t - x_{i}^t  \right).
\end{equation}
As before, \cite{Rahmaniani2020} show that the cuts~\eqref{BDD:ExtraStrengthenedCuts} are only useful for a fractional candidate solution,  where the cuts~\eqref{BDD:ExtraStrengthenedCuts} are stronger than both the standard optimality cuts~\eqref{OptimalityCutSingle} and the cuts~\eqref{BDD:StrengthenedCuts}.

However, in order to apply the Benders dual decomposition approach, we must repeatedly solve Lagrangian subproblems. Notably, depending on the structure of $\Omega$ and the sizes of $I$, solving the Lagrangian subproblems~\eqref{BDD:BendersSubproblemLagrangianPartial} or~\eqref{BDD:BendersSubproblemLagrangianFull} may be computationally infeasible. As such, viability of this method is dependent on the application.

\paragraph{Partial Benders Decomposition}
\label{SectionPartialBenders}

In Partial Benders Decomposition, information from the Benders subproblem is added to the main problem to improve the dual bound. This is done via including some variables of the Benders subproblem in the main problem. Following the work in \cite{Crainic2021}, we partition $J$ into two sets, $\bar{J}$ and $J\setminus \bar{J}$, and define the following reformulation of the maximum covering problem:
\begin{subequations}
\begin{alignat}{3}
\operatorname{Maximise } & \, \, \sum_{t = 1}^T \left( \sum_{j \in \bar{J}} d_j^t z_j^t + \sum_{j \in J\setminus \bar{J}} d_j^t z_j^t\right), 
\\
\text{subject to }
& \eqref{ModelMaximumCover:Domain}, \eqref{ModelMaximumCover:BinaryX}
\\ 
& \sum_{i \in I}  a_{ij}^{t} x_j^t \geq z_j^t, && 1 \leq t \leq T, j \in \bar{J},
\\
& \sum_{i \in I}  a_{ij}^{t} x_j^t \geq z_j^t, && 1 \leq t \leq T, j \in J \setminus \bar{J},
\\
& z_j^t \in \brace{0, 1} && 1 \leq t \leq T, j \in \bar{J},
\\
& z_j^t \in \brace{0, 1} && 1 \leq t \leq T, j \in J \setminus \bar{J}.
\end{alignat}
\label{PBD:ModelMaximumCover}
\end{subequations}
This formulation is, clearly, equivalent to the original formulation~\eqref{ModelMaximumCover}. However, when we project the $z$ variables for the Benders decomposition, we keep the variables $z_j^{t}, j \in \bar{J}$ in the main problem.  More specifically, in iteration $v$ of the Benders decomposition method, our main problem is given by 

\begin{subequations}
\begin{alignat}{3}
\operatorname{Maximise } & \, \, \theta + \sum_{t = 1}^T \sum_{j \in \bar{J}} d_j^t z_j^t, 
\\
\text{subject to }
& \eqref{ModelMaximumCover:Domain}, \eqref{ModelMaximumCover:BinaryX}, \eqref{SingleCut:BendersMain:OptimalityCuts}, \eqref{SingleCut:BendersMain:DomainTheta}
\\ 
& \sum_{t = 1}^T \sum_{i \in I}  a_{ij}^{t} x_j^t \geq z_j^t, && 1 \leq t \leq T, j \in \bar{J},
\label{PBD:BendersMainProblem:Covering}
\\
& z_j^t \in \brace{0, 1} && 1 \leq t \leq T, j \in \bar{J},
\end{alignat}
\label{PBD:BendersMainProblem}
\end{subequations}

and the Benders primal subproblem is then given by
\begin{subequations}
	\begin{alignat}{5}
\operatorname{Maximise} & \, \, \sum_{t = 1}^T \sum_{j \in J \setminus \bar{J}} d_j^t z_j^t,
\\ \text{subject to }
& z_j^t \leq \tilde{I}_j^t, && 1 \leq t \leq T, j \in J \setminus \bar{J}, 
\\
& 
z_j^t \in [0,1], && 1 \leq t \leq T, j \in J \setminus \bar{J},
\end{alignat}
\label{BendersSubproblemSingle}
\end{subequations}
for which optimality cuts can be derived as before. 

The users $j'$ in $J_s$ are particularly well-suited for scenario retention, i.e. to make $\bar{J} = J_s$. If we take $j' \in J_s$, by definition,, we have that $\sum_{t=1}^T \sum_{i \in I} a_{ij'}^t = 1$. Next, we take the period $t'$ and the facility $i'$ such that $a_{i'j'}^{t'} = 1$. Then, since $d_j^t >0$ by definition, we have that Constraint~\eqref{PBD:BendersMainProblem:Covering} will be satisfied with equality. As such, the variable $z_{j'}^{t'}$ can be removed from the main program~\eqref{PBD:BendersMainProblem}, and the objective function can be replaced by 
$$
\theta + \sum_{t = 1}^T \sum_{i \in I} \left( \sum_{j \in J_s} a_{ij}^t d_j^t \right) x_i^t.
$$

In addition to the scenario-retention method proposed in \cite{Crainic2021}, the authors also propose a scenario-creation method. This involves the creation of artificial scenarios given by a convex combination of real scenarios.  These artificial scenarios are then kept in the main problem as a proxy for their subproblem counterparts. However, this method relies on the variables $z_j^t$, which are no longer present in the branch-and-Benders-cut framework presented in Section~\ref{SectionSingleCut}. As such, this approach is not suited for this framework.

\section{Local Branching for Branch-And-Benders-Cut}
\label{SectionLocalBranching}

In the local branching method by \cite{Fischetti2003} for tackling mixed-integer programs, small subdomains of the feasible space are defined via distance-based neighbourhoods around solutions. These small subdomains are then solved in a separate subproblem via a black-box solver, and excluded from the feasible space, thus gradually reducing the size of the search space. In \cite{Rei2009}, this process was applied to a branch-and-Benders-cut framework, with the goal of simultaneously improving the upper and lower bounds for the search tree. 

We start in Section~\ref{SectionLocalBranchingOverview} by presenting the framework~\citep[following the work in ][]{Rei2009} for embedding a local branching approach within our branch-and-Benders-cut method. This defines a modified main problem and the necessary subproblems, which are all reliant on a distance metric. In Section~\ref{SectionDistanceMetric}, we propose a new distance metric for the dynamic MCLP. Via this distance metric, we then provide a novel solution method for quickly solving subproblems. Finally, in Section~\ref{SectionBranchingScheme}, we discuss methods for separating the feasible subdomains explored in our subproblems.

\subsection{Overview}
\label{SectionLocalBranchingOverview}

We implement the local branching scheme within the framework of our branch-and-Benders-cut methods, as proposed in \cite{Rei2009}. Since it is based on our branch-and-Benders-cut methods, any of the acceleration techniques presented in Section~\ref{SectionImprovements} can be applied to the local branching method as well. To that end, contrary to Section~\ref{SectionImprovements}, we present this section using the multi-cut version of the Benders optimality cuts. Due to the central role of the time period in the proceeding developments, the multi-cut formulation allows for a more natural and interpretable explanation. Thus, at iteration $v$ of the Benders decomposition, we consider the following local branching main problem:
\begin{subequations}
\begin{alignat}{3}
\operatorname{Maximise} \quad & \sum_{t=1}^T \theta^t, 
\label{LocalBranching:MainProblem:Objective} 
\\
\text{subject to} \quad 
&\eqref{ModelMaximumCover:Domain}, \eqref{ModelMaximumCover:BinaryX},\eqref{SingleCut:BendersMain:DomainTheta}
\notag
\\
& \sum_{i \in I} \left( \sum_{j \in \Gamma^t(x^{rt})} d_j^t a_{ij}^t \right) x_j^t \geq \theta^t - \sum_{j \in J\setminus \Gamma^t(x^{rt})} d_j^t, \quad && 1 \leq t \leq T, 1 \leq r \leq v,
\label{LocalBranching:MainProblem:OptimalityCuts}
\\
& \operatorname{Dist}_{x^s} \left( x \right) \geq {\kappa}^s, \quad && 1 \leq s \leq u,
\label{LocalBranching:MainProblem:Distance}
\end{alignat}
\label{LocalBranching:MainProblem}
\end{subequations}
where $x^s$ for $1\leq s \leq u$ are feasible facility location decisions previously determined. The function $\operatorname{Dist}_{x^s}$  represents a distance metric, indicating that only solutions found at a distance ${\kappa}^s$ from solution $x^s$ may be considered. We note that Constraints~\eqref{LocalBranching:MainProblem:OptimalityCuts} are the Benders optimality cuts~\eqref{MultiCut:BendersMain:OptimalityCuts} provided by Proposition~\ref{PropositionMulticut}.

Our goal, given an integer candidate solution $\tilde{x}$ and a threshold distance $\tilde{\kappa}$, is to find the optimal solution to our original problem~\eqref{ModelMaximumCover} restricted to distance $\tilde{\kappa}$ around $\tilde{x}$. As described in \cite{Rei2009}, we can then generate Benders optimality cuts~\eqref{OptimalityCutMulti} for this high-quality solution, thus improving the upper bound in the local branching main problem~\eqref{LocalBranching:MainProblem}. Simultaneously, the solution may also improve upon the incumbent, thus also increasing the lower bound. We can then exclude the subdomain of distance $\tilde{\kappa}$ around $\tilde{x}$ from the local branching main problem~\eqref{LocalBranching:MainProblem} by adding a constraint of type~\eqref{LocalBranching:MainProblem:Distance}.

Formally, let $\hat{x}$ denote the optimal solution to the \emph{restricted subproblem} centered around $\tilde{x}$, given by the solution to the following subproblem:
\begin{subequations}
\begin{alignat}{3}
    \operatorname{Maximise} \, \,& \sum_{t=1}^T \sum_{j \in J} d_j^t z_j^t, \\
    \text{subject to} \, \, &\eqref{ModelMaximumCover:Domain}, \eqref{ModelMaximumCover:BinaryX}, \eqref{SingleCut:BendersSubproblemPrimal:LowerBoundZ}, \eqref{SingleCut:BendersSubproblemPrimal:UpperBoundZ}, \eqref{LocalBranching:MainProblem:Distance} \notag,
    \\
    & \operatorname{Dist}_{\tilde{x}} \left( x \right) \leq \tilde{\kappa}.
    \label{LocalBranching:RestrictedSubproblem:Distance}
\end{alignat}
\label{LocalBranching:RestrictedSubproblem}
\end{subequations}

If $\hat{x}$ has a strictly higher objective than $\tilde{x}$, we generate  Benders optimality cuts~\eqref{OptimalityCutMulti} for $\tilde{x}$, and we set $x^{s+1} = \tilde{x}, \kappa^{s+1} = \tilde{\kappa} + 1, s = s+1$ and repeat the restricted problem~\eqref{LocalBranching:RestrictedSubproblem} centered around $\hat{x}$. If $\hat{x} = \tilde{x}$, we create a \emph{diversified subproblem} which replaces the threshold distance $\tilde{\kappa}$ in Constraint~\eqref{LocalBranching:RestrictedSubproblem:Distance} with $\tilde{\kappa}' > \tilde{\kappa}$, and adds the constraint $\operatorname{Dist}_{\tilde{x}} \left( x \right) \geq 1$. This guarantees that the resulting optimal solution $\hat{x}$ is different than $\tilde{x}$, and we create Benders optimality cuts~\eqref{OptimalityCutMulti} from $\hat{x}$. 

By setting $x^{s+1} = \tilde{x}, \kappa^{s+1} = \tilde{\kappa} + 1$, we ensure that the feasible domain in the restricted subproblem~\eqref{LocalBranching:RestrictedSubproblem} in iteration $v$ and the main problem in iteration $v+1$ are complementary. More specifically, the restricted subproblem only considers solutions with distance $\operatorname{Dist}_{\tilde{x}} \left( x \right) \leq \tilde{\kappa}$, whereas the main problem considers solutions with distance $\operatorname{Dist}_{\tilde{x}} \left( x \right) \geq \tilde{\kappa} + 1$. Since the candidate solution $\tilde{x}$ is integer feasible, it is not necessary to consider solutions with distance $\operatorname{Dist}_{\tilde{x}} \left( x \right) \in (\tilde{\kappa}, \tilde{\kappa}+1)$, as these correspond to fractional solutions.

We remark that the local branching main problem~\eqref{LocalBranching:MainProblem} and the restricted subproblem~\eqref{LocalBranching:RestrictedSubproblem} can be defined via any distance metric and with any threshold distance $\tilde{\kappa}$. By increasing the size of the subdomain, we can remove a larger area from the feasible space of the local branching main problem~\eqref{LocalBranching:MainProblem}, but at the cost of a restricted subproblem~\eqref{LocalBranching:RestrictedSubproblem}  that is harder to solve than that of a smaller size. We also note that any method can be used to find the optimal solution $\hat{x}$ to the restricted subproblem~\eqref{LocalBranching:RestrictedSubproblem}. In \cite{Fischetti2003}, a general-purpose mixed-integer linear solver is used, while the branch-and-Benders-cut approach with the acceleration techniques proposed in Section~\ref{SectionImprovements} is a natural choice in our case. However, by carefully selecting our distance metric and threshold distance $\tilde{\kappa}$, we can derive an exact solution method for the restriced subproblem~\eqref{LocalBranching:RestrictedSubproblem} which is more effective than the accelerated branch-and-Benders-cut method.

\subsection{An effective formulation for a tailored distance metric}
\label{SectionDistanceMetric}

The distance metric typically used is the Hamming distance~\citep[see, e.g., ][]{Fischetti2003, Rei2009}. This results in Constraint~\eqref{LocalBranching:RestrictedSubproblem:Distance} taking the form
\begin{equation}
    \sum_{t=1}^T \sum_{i \in I: \tilde{x}_i^t = 1} \left( 1 - x_i^t\right) + \sum_{t=1}^T\sum_{i \in I: \tilde{x}_i^{t} = 0} x_i^t \leq \tilde{\kappa}.
    \label{LocalBranching:HammingDistance}
\end{equation}

However, we propose to use a modified distance metric in the dynamic case, which enforces that the Hamming distance in each time period must be within our threshold. This results in Constraint~\eqref{LocalBranching:RestrictedSubproblem:Distance} taking the form
\begin{equation}
    \sum_{i \in I: \tilde{x}_i^t = 1} \left( 1 - x_i^t\right) + \sum_{i \in I: \tilde{x}_i^{t} = 0} x_i^t \leq \tilde{\kappa}, \, \, 1 \leq t \leq T.
    \label{LocalBranching:NewDistance}
\end{equation}

The motivation behind this new distance metric derives from infrastructure contexts, where the facilities under consideration correspond to significant investments (e.g. warehouses, stores, etc.). In those contexts, facilities which are added in early time periods are likely to persist throughout the time horizon. In the case of the general Hamming distance, this incurs a repeated penalty in each time period, increasing the total distance. As an illustration, consider the example in Figure~\ref{ExampleOverlap} and the two solutions $x^t = (1,0,0)$ and $\hat{x}^t = (1,1,0), 1 \leq t \leq T$, and a time horizon $T=4$. These solutions are at distance $4$ when considering the general Hamming distance, whereas they are only at distance $1$ using the new metric.

An important benefit of this new distance metric allows for an efficient solving method if we consider $\tilde{\kappa} = 2$. To describe this method we first note that, for a candidate solution $\tilde{x}$ and time period $t$, there are only two sets of modifications which result in an integer solution at distance exactly $1$:
\begin{enumerate}[label=1\alph*)]
    \item Add one facility which is currently not selected. \label{EnumD1Add}
    \item Remove one facility which is currently selected. \label{EnumD1Remove}
\end{enumerate}
Likewise, at distance exactly $2$, there are only three sets of modifications:
\begin{enumerate}[label=2\alph*)]
    \item Add two facilities which are currently not selected. \label{EnumD2Add}
    \item Remove one facility which is currently selected, and add one facility which is currently not selected. 
    \label{EnumD2Swap}
    \item Remove two facilities which are currently selected.
    \label{EnumD2Remove}
\end{enumerate}
Since, by assumption $a_{ij}^t \in \brace{0,1}$ and $d_j^t > 0$, options~\ref{EnumD1Remove} and~\ref{EnumD2Remove} cannot lead to an increase in the objective function, and hence, can be disregarded. 

As a consequence, there is a limited number of feasible solutions contained within the restricted subproblem~\eqref{LocalBranching:RestrictedSubproblem} when considering $\tilde{\kappa} = 2$. We can then use Proposition~\ref{TheoremRestrictedSubproblemSolution} to evaluate the quality of feasible solutions (as given in Equation~\eqref{SolutionCoverage}) based on Benders optimality cuts.

\begin{boxedproposition}
\label{TheoremRestrictedSubproblemSolution}
Let $x$ be an (integer) feasible solution to the maximum covering model~\eqref{ModelMaximumCover}, and let $t$ be any time period, $1 \leq t \leq T$. Let $\hat{i} \in I$ be such that $x_{\hat{i}}^{t} = 0$. Then the modified solution $\hat{x}$ with $\hat{x}_{\hat{i}}^t =1$ and $\hat{x}_i^t = x_i^t, i \neq \hat{i}$ satisfies

$$
\sum_{i \in I} \left( \sum_{j \in \Gamma^t(x^t)} d_j^t a_{ij}^t\right) \hat{x}_i^t + \sum_{j \in J \setminus \Gamma^t(x^t)} d_j^t = \sum_{j \in J} \min \left\{ 1, \sum_{i \in I} a_{ij}^t \hat{x}_i^t\right\},
$$
when $\Gamma^t(x^t)$ is given by \eqref{SingleB0} or \eqref{SingleB1}. 
\end{boxedproposition}
\begin{proof}
We first note that, by rearranging the terms in the left-hand side, we have that 
\begin{equation*}
        \sum_{i \in I} \left( \sum_{j \in \Gamma^{t}(x)} d_j^{t} a_{ij}^{t}\right) \hat{x}_i^{t} + \sum_{j \in J \setminus \Gamma^{t}(x)} d_j^{t}  = \sum_{j \in \Gamma^{t}(x)} \left( \sum_{i \in I} a_{ij}^{t} \hat{x}_i^t\right) d_j^{t} +  \sum_{j \in J \setminus \Gamma^{t}(x)} d_j^{t}.
\end{equation*}
Now, we concentrate on both cases for $\Gamma^{t}(x)$.

If $\Gamma^{t}(x)$ is given by \eqref{SingleB0}, then we have that
\begin{align*}
    \sum_{j \in \Gamma^{t}(x)} \left( \sum_{i \in I} a_{ij}^{t} \hat{x}_i^t\right) d_j^{t} +  \sum_{j \in J \setminus \Gamma^{t}(x)} d_j^{t} & = \sum_{j \in J: I(x)_j^{t} < 1} \left( \sum_{i \in I } a_{ij}^{t} \hat{x}_i^{t}\right) d_j^{t} + \sum_{j \in J: I(x)_j^{t} \geq 1} d_j^{t}, \\
    & = \sum_{j \in J: I(x)_j^{t} = 0} \left( \sum_{i \in I\setminus \{i\}} a_{ij}^{t} x_i^{t}  +  a_{\hat{i}j}^{t} \right) d_j^{t} + \sum_{j \in J: I(x)_j^{t} \geq 1} d_j^{t},
\end{align*}
where $I(x)_j^{t} = \sum_{i \in I} a_{ij}^{t} x_i^{t} < 1$ implies that $I(x)_j^{t} =0$, since all elements in the summation are binary.

From there, for the first term, $a_{\hat{i}j}^{t}$ being binary also implies that $ \sum_{i \in I} a_{ij}^{t} \hat{x}_i^t = \sum_{i \in I\setminus \{\hat{i}\}} a_{ij}^{t} x_i^{t}  +  a_{\hat{i}j}^{t} = I_j^t(x)  +  a_{\hat{i}j}^{t}\leq 1$. As such, we have that 
$$ \sum_{j \in J: I(x)_j^{t} = 0} \left( \sum_{i \in I\setminus \{\hat{i}\}} a_{ij}^{t} x_i^{t}  +  a_{\hat{i}j}^{t} \right) d_j^{t} = \sum_{j \in J: I(x)_j^{t} = 0} \min \{1, \sum_{i \in I} a_{ij}^{t} \hat{x}_i^t \}d_j^{t}$$.
For the second term, we note that since $\hat{x}$ has an extra facility compared to $x$, we also have that the coverage for each user must be equal or greater. Concretely, $I(x)_j^{t} = \sum_{i \in I} a_{ij}^{t} x_i^t \leq \sum_{i \in I} a_{ij}^{t} \hat{x}_i^t$, and thus $\sum_{j \in J: I(x)_j^{t} \geq 1} d_j^{t} = \sum_{j \in J: I(x)_j^{t} \geq 1} \min \{1, \sum_{i \in I} a_{ij}^{t} \hat{x}_i^t \} d_j^{t}$. 

As a consequence, we have that
\begin{align*}
    \sum_{j \in J: I(x)_j^{t} = 0} \left( \sum_{i \in I\setminus \{\hat{i}\}} a_{ij}^{t} x_i^{t}  +  a_{\hat{i}j}^{t} \right) d_j^{t} + \sum_{j \in J: I(x)_j^{t} \geq 1} d_j^{t} & = \sum_{j \in J} \min \{1, \sum_{i \in I} a_{ij}^{t} \hat{x}_i^t \}d_j^{t}.
\end{align*}

If $\Gamma^{t}(x)$ is given by \eqref{SingleB1}, then we have that
\begin{align*}
    & \sum_{j \in \Gamma^{t}(x)} \left( \sum_{i \in I} a_{ij}^{t} \hat{x}_i^t\right) d_j^{t} +  \sum_{j \in J \setminus \Gamma^{t}(x)} d_j^{t} 
    \\
   = & \sum_{j \in J_s} \left( \sum_{i \in I} a_{ij}^{t} \hat{x}_i^t \right)d_j^t +  \sum_{j \in J\setminus J_s: I(x)_j^{t} < 1} \left( \sum_{i \in I } a_{ij}^{t} \hat{x}_i^{t}\right) d_j^{t} + \sum_{j \in J\setminus J_s: I(x)_j^{t} \geq 1} d_j^{t}.
\end{align*}
We note that, by definition of the set $J_s$, we always have that $\sum_{i \in I} a_{ij}^{t} \hat{x}_i^t \leq 1$. As a consequence, we have that $\sum_{j \in J_s} \left(\sum_{i \in I} a_{ij}^{t} \hat{x}_i^t \right)d_j^t = \sum_{j \in J_s} \min \{ 1, \sum_{i \in I} a_{ij}^{t} \hat{x}_i^t \}d_j^t$. The last two terms can be transformed using reasoning similar to that of the prior case, which gives our result. 

\end{proof}

We remark that Proposition~\ref{TheoremRestrictedSubproblemSolution} does not hold in general if $\Gamma^t(x^t)$ is given by \eqref{SingleB2}. As a counterexample if $\Gamma^{t}(x^t)$ is given by \eqref{SingleB2}, consider the example illustrated in Figure~\ref{ExampleOverlap}, the candidate solution $x = (1,0,0)$, and $\hat{i} = i_2$. We then have that $\sum_{i \in I} \left( \sum_{j \in \Gamma^{t}(x^t)} d_j^{t} a_{ij}^{t}\right) \hat{x}_j^{t} + \sum_{j \in J \setminus \Gamma^{t}(x^t)} d_j^{t} = 36$ and $\sum_{t=1}^T \sum_{j \in J} \min \left\{ 1, \sum_{i \in I} a_{ij}^t \hat{x}_i^t\right\} = 28$.

Put more simply, Proposition~\ref{TheoremRestrictedSubproblemSolution} states that if we take a given solution $x$ and only add one facility to it, then the Benders optimality cuts for $x$ accurately give the coverage for the modified solution. 

We can then derive an efficient method for finding the optimal solution of distance two around $\tilde{x}$ by combining the limited number of valid modifications to $\tilde{x}$ along with Proposition~\ref{TheoremRestrictedSubproblemSolution}. For this, let $e^{i}$ denote the elemental vector with $1$ in position $i$ and $0$ elsewhere. We then denote $x^t + e^i$ (respectively, $x^t - e^i$) as the modified solution which adds the currently unused facility $i$ in period $t$ (respectively, removes the currently used facility $i$).  We present this method in 
Proposition~\eqref{TheoremSubproblemMethod}.

\begin{boxedproposition}
\label{TheoremSubproblemMethod}
Let $\hat{x}_j^t$ denote the optimal values for the variables $x_j^t$ in the restricted subproblem~\eqref{LocalBranching:RestrictedSubproblem} around an integer candidate solution $\tilde{x}$ at distance $\tilde{\kappa} = 2$. Then  $\hat{x}_j^t$ are also the optimal values for the variables $x_j^t$ in the following problem:
\begin{subequations}
	\begin{alignat}{5}
\operatorname{Maximise } & \, \, \sum_{t=1}^T \theta^t, 
\\
\operatorname{subject to }
& \eqref{ModelMaximumCover:Domain}, \eqref{ModelMaximumCover:BinaryX}, \eqref{SingleCut:BendersMain:OptimalityCuts},  \eqref{SingleCut:BendersMain:DomainTheta}, \eqref{LocalBranching:MainProblem:Distance}, \notag
\\ 
& \sum_{i \in I: \tilde{x}_i^t = 1} \left( 1 - x_i^t\right) + \sum_{i \in I: \tilde{x}_i^{t} = 0} x_i^t \leq 2, && 1 \, \leq t \leq T,
\label{LocalBranching:RestrictedSubproblemReformulated:Distance}
\\
&\sum_{i \in I} \left( \sum_{j \in J\setminus \Gamma^t (\tilde{x}^t)} a_{ij}^{t} d_j^t \right) x_{i}^t + \sum_{j \in \Gamma^t (\tilde{x}^t)} d_j^t \geq \theta^t, \, && 1\leq t \leq T, 
\label{LocalBranching:RestrictedSubproblemReformulated:TrustAdd1Cuts}
\\
&\sum_{i \in I} \left( \sum_{j \in J\setminus \Gamma^t (\tilde{x}^t + e^{\hat{i}}) } a_{ij}^{t} d_j^t \right) x_{i}^t + \sum_{j \in \Gamma^t (\tilde{x}^t + e^{\hat{i}})} d_j^t \geq \theta^t, \, && 1 \leq t \leq T, \hat{i} : \tilde{x}_{\hat{i}}^t = 0,
\label{LocalBranching:RestrictedSubproblemReformulated:TrustAdd2Cuts}
\\
&\sum_{i \in I} \left( \sum_{j \in J\setminus \Gamma^t (\tilde{x}^t - e^{\hat{i}})} a_{ij}^{t} d_j^t \right) x_{i}^t + \sum_{j \in \Gamma^t (\tilde{x}^t - e^{\hat{i}})} d_j^t \geq \theta^t, \, && 1 \leq t \leq T, \hat{i} : \tilde{x}_{\hat{i}}^t = 1.
\label{LocalBranching:RestrictedSubproblemReformulated:TrustSwapCuts} 
\end{alignat}
\label{LocalBranching:RestrictedSubproblemReformulated}
\end{subequations}
\end{boxedproposition}
\begin{proof}
The feasible space for the variables $x_j^t$ are defined via Constraints~\eqref{ModelMaximumCover:Domain},~\eqref{ModelMaximumCover:BinaryX},~\eqref{LocalBranching:MainProblem:Distance}, and~\eqref{LocalBranching:RestrictedSubproblemReformulated:Distance}. Since these constraints are present in both models, it follows that the set of feasible values for $x_j^t$ are the same. We now show that the objective value for every feasible solution $x$ is the same in both models.

From our previous discussion, there is a limited number of feasible solutions which we must consider: the solution $\tilde{x}$ itself, as well as the resulting solutions from modifications~\ref{EnumD1Add},~\ref{EnumD2Add} and~\ref{EnumD2Swap}. The value of the solution $\tilde{x}$ is given by Constraint~\eqref{LocalBranching:RestrictedSubproblemReformulated:TrustAdd1Cuts}. For all other solutions, we make use of Proposition~\ref{TheoremRestrictedSubproblemSolution}. More specifically, for each type of modification, we have generated an optimality cut for a solution with exactly one fewer facilities:
\begin{itemize} 
\item Constraints~\eqref{LocalBranching:RestrictedSubproblemReformulated:TrustAdd1Cuts} generate optimality cuts for the solution $\tilde{x}^t$ for all time periods $t$. As a consequence, all solutions which add one facility to $\tilde{x}$ (i.e. modifications of type~\ref{EnumD1Add} are evaluated correctly.
\item Constraints~\eqref{LocalBranching:RestrictedSubproblemReformulated:TrustAdd2Cuts} generate optimality cuts for the solutions $\tilde{x}^t + e^{\hat{i}}$ for all time periods $t$ and for all facilities $\hat{i}$ which are currently not installed. As a consequence, all solutions which add one additional facility beyond $\tilde{x}^t + e^{\hat{i}}$ (thus modifications of type~\ref{EnumD2Add} are evaluated correctly.
\item Constraints~\eqref{LocalBranching:RestrictedSubproblemReformulated:TrustSwapCuts}  generate optimality cuts for the solutions $\tilde{x}^t - e^{\hat{i}}$ for all time periods $t$ and for all facilities $\hat{i}$ which are currently not installed. As a consequence, all solutions removing $\hat{i}$ and then adding another facility (thus modifications of type~\ref{EnumD2Swap} are evaluated correctly.
\end{itemize}

\end{proof}

 We note the following observations about Proposition~\ref{TheoremSubproblemMethod}:
 \begin{itemize}
      \item Due to the presence of Constraints~\eqref{ModelMaximumCover:BinaryX}, model~\eqref{LocalBranching:RestrictedSubproblemReformulated} is a mixed-integer linear program, which can be solved directly via a generic MILP solver. When solving model~\eqref{LocalBranching:RestrictedSubproblemReformulated} directly with a generic solver is the strategy for finding the optimal solution to the restricted subproblem~\eqref{LocalBranching:RestrictedSubproblem}, we denote this as \subD. When, instead, the restricted subproblem~\eqref{LocalBranching:RestrictedSubproblem} is solved via the branch-and-Benders-cut method, we denote this \subB. 
      
     \item Some of Constraints~\eqref{LocalBranching:RestrictedSubproblemReformulated:TrustAdd1Cuts}-\eqref{LocalBranching:RestrictedSubproblemReformulated:TrustSwapCuts} may be redundant, in the sense that they allow for the evaluation of solutions which can either be evaluated using other constraints within the model or which correspond to infeasible solutions. 
     In general, verifying the feasibility of the solutions before generating the associated Constraints~ \eqref{LocalBranching:RestrictedSubproblemReformulated:TrustAdd2Cuts} or~\eqref{LocalBranching:RestrictedSubproblemReformulated:TrustSwapCuts} is unlikely to be beneficial, due to the speed at which these constraints can be built. Instead, infeasible solutions can be detected by the MILP solver via  Constraints~\eqref{ModelMaximumCover:Domain} and~\eqref{LocalBranching:MainProblem:Distance}.  However, some structures of $\Omega$ may make the detection of some infeasible solutions a trivial task, such as the use of precedence constraints, and thus justify the computational effort of the verification process. 
     
     \item From the proof of Proposition~\ref{TheoremSubproblemMethod}, it is clear that the optimal objective value of the restricted subproblem~\eqref{LocalBranching:RestrictedSubproblem} and of model~\eqref{LocalBranching:RestrictedSubproblemReformulated} are the same. This is important, as it allows us to verify if the newly found solution has a better objective value than the incumbent, and update the incumbent accordingly.
     
 \end{itemize}

\subsection{Branching}
\label{SectionBranchingScheme}
Using Proposition~\ref{TheoremSubproblemMethod}, we can thus easily solve the restricted subproblem~\eqref{LocalBranching:RestrictedSubproblem}. However, this ease of solving comes at a cost when moving to solve the full problem~\eqref{LocalBranching:MainProblem}.

After each candidate solution $\tilde{x}$ and the set of restricted and diversified subproblems, the Constraints~\eqref{LocalBranching:MainProblem:Distance} ensure that the solver removes the subdomain examined during the subproblems, thus avoiding unnecessary work. In \cite{Rei2009}, the Hamming distance is used in Constraints~\eqref{LocalBranching:MainProblem:Distance}, and hence, the complement of the set $\tilde{\Omega} = \brace{x \in \Omega: \eqref{ModelMaximumCover:BinaryX}, \operatorname{Dist}_{\tilde{x}} \left( x \right) \leq \tilde{\kappa}}$ is simply 
$$\tilde{\Omega}^{\mathbf{c}} = \brace{x \in \Omega: \eqref{ModelMaximumCover:BinaryX}, \sum_{t=1}^T \sum_{i \in I: \tilde{x}_i^t = 1} \left( 1 - x_i^t\right) + \sum_{i \in I: \tilde{x}_i^{t} = 0} x_i^t \geq \tilde{\kappa} +1}.$$

However, using the distance metric in Section~\ref{SectionDistanceMetric}, the set $\tilde{\Omega}^{\mathbf{c}}$ requires the Hamming distance to be at least $\tilde{\kappa} +1$ in at least one time period, which cannot be modeled as a single linear constraint. Therefore, we must create a branch for each time period in order to mimic $\tilde{\Omega}^{\mathbf{c}}$ through $T$ disjoint problems. More specifically, for every $1 \leq t' \leq T$, we generate the following problem:
\begin{subequations}
\begin{alignat}{3}
\operatorname{Maximise} \quad & \sum_{t=1}^T \theta^t, 
\label{LocalBranching:Branching:Branches:Objective} 
\\
\text{subject to} \quad 
&\eqref{ModelMaximumCover:Domain}, \eqref{ModelMaximumCover:BinaryX},\eqref{SingleCut:BendersMain:DomainTheta}, \eqref{LocalBranching:MainProblem:OptimalityCuts}, \eqref{LocalBranching:MainProblem:Distance},
\notag
\\
& \sum_{i \in I: \tilde{x}_i^{t'} = 1} \left( 1 - x_i^{t'}\right) + \sum_{i \in I: \tilde{x}_i^{t'} = 0} x_i^{t'} \geq \tilde{\kappa} + 1, 
\label{LocalBranching:Branching:Branches:CurrentYear}
\\
& \sum_{i \in I: \tilde{x}_i^{t} = 1} \left( 1 - x_i^{t}\right) + \sum_{i \in I: \tilde{x}_i^{t} = 0} x_i^{t'} \leq \tilde{\kappa}, t < t'.
\label{LocalBranching:Branching:Branches:PriorYears}
\end{alignat}
\label{LocalBranching:Branching:Branches}
\end{subequations}
Each of these problems can be viewed as a branch at node $\tilde{x}$ in the search tree. In this way, the Constraints~\eqref{LocalBranching:MainProblem:Distance} are formed of the branching constraints~\eqref{LocalBranching:Branching:Branches:CurrentYear} and~\eqref{LocalBranching:Branching:Branches:PriorYears} from the prior nodes in the search tree.  In essence, these branches progressively select each time period as the one that exceeds the threshold distance, and impose all prior years to be below the threshold in order to ensure distinct sets. When using these branches to separate the subdomains of the restricted and diversified subproblems, we refer to the method as \sepB.

Alternatively, we can model this disjunction through a set of cuts, requiring the introduction of auxiliary binary variables and Big-M constraints. More specifically, the set $\tilde{\Omega}^{\mathbf{c}}$ can be obtained by adding to Constraints~\eqref{LocalBranching:MainProblem:Distance} the constraints
\begin{subequations}
\begin{alignat}{4}
    \sum_{i: \tilde{x}_{i}^t = 1 } (1 - x_{i}^t) +  \sum_{i: \tilde{x}_{i}^t = 0 } x_{i}^t + \tilde{\delta}^t (\tilde{\kappa} + 1) &\geq \tilde{\kappa} + 1, && \quad 1 \leq t \leq T, 
    \label{LocalBranching:DistanceOuter:BigM}
    \\
    \sum_{t=1}^T \tilde{\delta}^t &\leq T-1,
    \\
    \tilde{\delta}^t &\in \brace{0,1}, && \quad 1 \leq t \leq T.
\end{alignat}
\label{LocalBranching:DistanceOuter}
\end{subequations}
However, since we are adding the binary variables $\tilde{\delta}^t$ and the Big-M constraints~\eqref{LocalBranching:DistanceOuter:BigM} to the main program~\eqref{LocalBranching:MainProblem}, the problem gets progressively more difficult to solve as more subdomains are removed.

We note that in the case that the distance $\tilde{\kappa} = 0$, the Constraints~\eqref{LocalBranching:DistanceOuter} are equivalent to (and thus can be replaced by) the typical no-good cut:
\begin{equation}
    \sum_{t=1}^T  \sum_{i: \tilde{x}_{i}^t = 1 } (1 - x_{i}^t) +  \sum_{i: \tilde{x}_{i}^t = 0 } x_{i}^t \geq 1,
    \label{LocalBranching:DistanceNoGood}
\end{equation}
which does not require the use of binary auxiliary variables.

\section{Computational experiments}
\label{SectionResults}
In this section, we compare the performance of the variations of the Benders decomposition techniques presented in this work. We start in Section~\ref{SectionResultsApplication} by giving the background information for the problem instances, which are obtained from an electric vehicle (EV) charging station location model, and based on real-life data. In Section~\ref{SectionResultsBaseline}, we then examine the performance metrics when using our Benders decomposition methods to solve these problem instances, which clearly demonstrates the capabilities and limitations of each method on a practical example.  In Section~\ref{SectionResultsSubproblem}, we evaluate the options for the local branching method presented in Section~\ref{SectionDistanceMetric} for the restricted subproblem solution method. These experiments demonstrate the efficiency of solving our reformulated restricted subproblem~\eqref{LocalBranching:RestrictedSubproblemReformulated}, validating its use in Section~\ref{SectionResultsBaseline}. Similarly, in Section~\ref{SectionResultsSpace}, we compare the options for eliminating feasible space already explored within the local branching method, as presented in Section~\ref{SectionBranchingScheme}. Additionally, we discuss the technical limitations that face the implementation of the \sepB method, and the consequences of using either the \sepD or \sepB methods.

To more easily distinguish between the methodologies discussed in this paper, we present the following nomenclature:
\begin{enumerate}
    \item The unaccelerated branch-and-Benders-cut procedure in \cite{Cordeau2019} as described in Section~\ref{SectionSingleCut} is referred to as the \singlecut method. This procedure is unaccelerated in the sense that it does not include any of the acceleration techniques presented in Section~\ref{SectionImprovements}. In particular, this method generates Benders optimality cuts~\eqref{OptimalityCutSingle} at every integer and fractional candidate solution.
    
    \item The accelerated branch-and-Benders-cut procedure presented in Section~\ref{SectionImprovements} is referred to as the \multicut method. As some techniques are dependent on the instance characteristics, we describe the set of acceleration techniques in more detail in Section~\ref{SectionResultsApplication}. Notably, for the \multicut method, Benders optimality cuts are generated for every integer candidate solution, while Benders optimality cuts are generated for fractional candidate solutions only at the root node. 
    
    \item The accelerated branch-and-Benders-cut procedure with local branching presented in Section~\ref{SectionLocalBranching} is referred to as the \localbranching{Sub}{Sep} method, where \textsc{Sub} indicates the restricted subproblem solution method as described in Section~\ref{SectionDistanceMetric} and \textsc{Sep} indicates the subdomain separation scheme as described in Section~\ref{SectionBranchingScheme}. In all cases, we impose a time limit of 60 seconds for the solving of each restricted and diversified subproblem, after which the the best incumbent solution is used for the subsequent procedures.
    
    Regardless of the subproblem solution method and the subdomain separation scheme, almost all acceleration techniques are identical between the \multicut and \localbranching{Sub}{Seb} methods. The one exception is the generation of optimality cuts for fractional solutions. It was observed that, in the hard instances, both the \singlecut and \multicut methods spent a considerable amount of time in the root node generating optimality cuts for fractional solutions. Since the \textsc{SepB} separation scheme relies on branching, we push the solver to exit the root node more quickly. To do this, we solve the continuous relaxation of the model using the multi-cut method, then keep the resulting (fractional) optimality cuts in the formulation. We then solve the MILP formulation of the problem and only generate optimality cuts for integer solutions. This procedure is similar to the one proposed in \cite{Fortz2009}, where feasibility cuts for fractional solutions are only generated as part of solving the LP relaxation. 
\end{enumerate}

The tests were run on a server running Linux version 3.10, with an Intel Core i7-4790 CPU with eight virtual cores and 32 GB of RAM. The code is written in C++, and is publicly available\footnote{\url{https://github.com/StevenLamontagne/BendersDecomposition}}. We use CPLEX version 22.1.1 limited to a single thread, and with the Benders cuts implemented via the generic callback feature. Other than the removal of incompatible preprocessing and reformulation options in the cases which use the callback feature, the only notable parameters which are not set at default value are the memory management and numerical precision. The Eigen library~\citep{Eigen} is used for efficient dense and sparse matrix calculations, notably for the calculation of $\tilde{I}_j^t$ and the subsequent generation of cuts. In all tests, we impose a two-hour (7,200 second) time limit for solving each instance.

\subsection{Application: Electric vehicle charging station placement}
\label{SectionResultsApplication}

We apply the proposed methods to the instances of the electric vehicle (EV) charging station placement problem by \cite{Lamontagne2022}. We chose these problem instances because of their practical interest and the large number of users, which leads to large-scale dynamic MCLPs and therefore allows us to demonstrate the capabilities of our methodological contributions.

In this problem, there are two groups making sequential decisions over a multi-period time span. First, a decision maker which decides where to place public charging infrastructure. And second, users purchasing a vehicle in the given time period, who must elect between an EV or a conventional vehicle (a choice which depends on the charging network). If a user is covered by a charging station (in the maximum covering sense), then they elect to purchase an EV. The objective of the decision maker is to place the charging infrastructure in such a way as to maximise EV adoption.

In addition to determining which charging stations to open, this problem also aims to find the optimal sizing (i.e. number of charging outlets) of each open station. For this purpose, the decision variables are binary, indicating if a charging station $i$ has \emph{at least} $k$ outlets (with $k$ ranging from $1$ to an upper bound $m_i$). The set $\Omega$ is then composed of three sets of constraints: precedence constraints imposing for each station that in order to have at least $k+1$ outlets, one must first have at least $k$ outlets, constraints that forbid the removal of outlets between time periods, and budget constraints which force the cost of installing new outlets to be below a budget in each time period.

The test instances are divided into five datasets, each with different characteristics and each containing 20 instances. We refer to \cite{Lamontagne2022} for a detailed description. However we note that the Simple, Distance, and HomeCharging datasets are easier to solve than the LongSpan and Price datasets. Notably, no instances in the latter two datasets were solved exactly in \cite{Lamontagne2022}, with the best incumbent objective value in the branch-and-cut method (via CPLEX 12.10) being below that of the greedy method described in that paper. As a consequence, we refer to the set of instances in the Simple, Distance, and HomeCharging datasets as the ``easy instances'', while those from the LongSpan and Price datasets are called the ``hard instances''.

Next, we discuss the customisation of the acceleration techniques described in Section~\ref{SectionImprovements} used for this problem.
\begin{itemize}

    \item There are two preprocessing techniques we can use to eliminate users from consideration, without affecting the objective value. The first technique, as proposed in \cite{Legault2022}, is to eliminate any users for which $a_{ij}^t = 0, \forall j \in J$. Note that such users they can never be covered regardless of our efforts, implying that their elimination does not affect the objective. The second technique is to eliminate any user which is ``precovered'' (in the sense that there is an existing facility or option which guarantees coverage, regardless of decisions in the model). In particular, in the HomeCharging dataset, there are users who have access to a home charging system for their EV. As a consequence, a subset of those users will purchase an EV, regardless of the state of the public charging network, and can be eliminated from consideration.
    
    \item For a heuristic warmstart, we use the greedy method proposed in \cite{Lamontagne2022}. This is achieved via the warmstart feature present in CPLEX.

    \item We keep the users $j \in J_s$ within the main problem, as described in Section~\ref{SectionPartialBenders}. 
        
    \item For the Benders optimality cuts, we use the multicut equivalent of the Pareto-optimal cuts in Proposition~\ref{TheoremParetoOptimal}, using B1-type cuts. In other words, for $1 \leq t \leq T$, the sets $\Gamma^t(\tilde{x}^t)$ in Proposition~\ref{PropositionMulticut} are given by 
    \begin{equation}
    \label{OptimalityCutParetoMulti}
    \Gamma^t(\tilde{x}^t) = \left\{ j \in J\setminus J_s: \tilde{I}_j^t < 1 \right\} \cup \left\{ j \in J\setminus J_s: \tilde{I}_j^t = 1, \overset{c}{I}_j^t < 1 \right\}.
    \end{equation}
    The users $j \in J_s$ do not need to be considered, due to the partial Benders decomposition strategy mentioned above. We use the iterative method proposed in \cite{Papadakos2008} for estimating a core point.
    
    \item In the cases of the LongSpan and Price datasets, CPLEX is unable to solve even the linear programming relaxation of model~\eqref{ModelMaximumCover} as shown in \cite{Lamontagne2022}. As such, we do not use Benders dual decomposition, as repeatedly solving the Lagrangian subproblems~\eqref{BDD:BendersSubproblemLagrangianPartial} or~\eqref{BDD:BendersSubproblemLagrangianFull} would be computationally infeasible.

\end{itemize}

\subsection{Comparison of methodologies}
\label{SectionResultsBaseline}

In this section, we compare the performance of the \singlecut, \multicut, and \localbranching{SubD}{SepD} methods for solving the instances. As will be discussed later, these options for the local branching method are the best suited for our application.  Additionally, we compare these with a standard branch-and-cut tree via CPLEX, and the greedy method from \cite{Lamontagne2022}. These last two methods are referred to as the \branchandcut and \greedy methods, respectively. We note that the \multicut and \localbranching{SubD}{SepD} methods both use the \greedy method as a warmstart, and thus will always terminate with an objective value that is at least equal to it.  

In Table~\ref{TableBaselineSummary}, we summarise our computational results for each methodology and instance set. The solving time, objective value, optimality gap, and number of nodes are reported as given by CPLEX. However, we note that the number of nodes in the local branching case only includes the ones for the main problem.  By examining the results, we note that:
\begin{itemize}
    
    \item While the \singlecut method is able to find a non-trivial solution (i.e. a solution which places at least one facility) for all instances, it found worse quality solutions on average than the \greedy method, \multicut, and \localbranching{SubD}{SepD} methods in the hard instances.
    
    \item The \multicut method was able to outperform the \singlecut method in all instances. In addition, in the easy instances, the \multicut method is able to run faster than the \branchandcut method. However, while it is able to reduce the optimality gap in the hard instances compared to the \singlecut method, the objective value is only slightly better than the \greedy method. In fact, the \multicut method only improved the objective value compared to the \greedy method in two instances in the LongSpan dataset, and two instances in the Price dataset.
    
    \item The \localbranching{SubD}{SepD} method required more time than either the \multicut or the \branchandcut methods in the easy instances, and the optimality gap is considerably higher than the \multicut method in the hard instances.  However, the \localbranching{SubD}{SepD} method improves the objective value compared to the \greedy method in all 20 instances in the LongSpan dataset and 17 instances in the Price dataset, which is notably more than the \multicut method.
\end{itemize}

Additional performance details are presented in \ref{AppendixResultsBaseline}.
\begin{table}
\centering
\resizebox{0.9\textwidth}{!}{\begin{tabular}{llrrrrr}
\hline
                     &                  &    Simple &  Distance &  HomeCharging &     LongSpan &        Price \\
\hline
Solve time (sec) & \greedy &      $<0.01$ &      $<0.01$ &          $<0.01$ &         0.05 &         0.10 \\
                     & \branchandcut &      0.20 &      0.70 &          6.56 &      7186.12 &      7194.97 \\
                     & \singlecut &      0.15 &      1.57 &         68.24 &      7184.25 &      7184.76 \\ 
                     & \multicut &      0.09 &      0.52 &          1.34 &      7184.54 &      7185.04 \\
                     & \localbranching{SubD}{SepD} &      4.67 &      4.81 &         13.00 &      7291.34 &      7305.71 \\

\hline 
Objective value & \greedy &  31814.20 &  16591.75 &      18016.47 &    133724.37 &     33641.74 \\
                     & \branchandcut &  31820.15 &  16627.14 &      18030.50 &         0.00 &         0.00 \\
                     & \singlecut &  31820.15 &  16627.14 &      18030.47 &    116663.50 &     32563.31 \\                      
                     & \multicut &  31820.15 &  16627.14 &      18030.40 &    133728.65 &     33642.19 \\
                     & \localbranching{SubD}{SepD} &  31820.15 &  16627.14 &      18030.49 &    133781.88 &     33650.81 \\
\hline 
Optimality gap (\%) & \branchandcut &      $<0.01$ &      $<0.01$ &          $<0.01$ & - & - \\
                     & \singlecut &      $<0.01$ &      $<0.01$ &          $<0.01$ &  15.79 &        15.00 \\ 
                     & \multicut &      $<0.01$ &      $<0.01$ &          $<0.01$ &         6.12 &        11.32 \\
                     & \localbranching{SubD}{SepD} &      $<0.01$ &      $<0.01$ &          $<0.01$ &        11.86 &        18.66 \\
\hline 
Number of nodes & \branchandcut &     59.25 &     43.35 &        131.95 &         0.00 &         0.00 \\
                     & \singlecut  &     52.40 &     27.85 &        228.50 &         0.00 &         0.00 \\ 
                     & \multicut &     65.35 &     75.05 &        860.45 &         0.00 &      1499.60 \\
                     & \localbranching{SubD}{SepD} &     65.05 &    205.60 &        797.35 &      7282.65 &      1202.50 \\
\hline 
\end{tabular}}
\caption{Average performance details for solution methods.}
\label{TableBaselineSummary}
\end{table}

\subsection{Further results: Local branching restricted subproblem solution method}
\label{SectionResultsSubproblem}

In this section, we compare the methods for solving the restricted subproblems~\eqref{LocalBranching:RestrictedSubproblem} of the local branching procedure. Since these subproblems are solved many times in each instance, it is imperative that they are solved efficiently. For this, we use the \subdirect and \subbenders procedures presented in Section~\ref{SectionDistanceMetric}. We recall that the \subdirect method solves the restricted subproblem~\eqref{LocalBranching:RestrictedSubproblem} by directly providing the reformulation~\eqref{LocalBranching:RestrictedSubproblemReformulated} to a generic MILP solver (in this case, CPLEX), while the \subbenders procedure solves the restricted subproblem~\eqref{LocalBranching:RestrictedSubproblem} with the \multicut method.

For a more accurate comparison between the \subdirect and \subbenders methods, we simulate the entire solving process. This starts with the MP~\eqref{LocalBranching:MainProblem}, with no optimality cuts~\eqref{MultiCut:BendersMain:OptimalityCuts} nor feasible space reductions~\eqref{LocalBranching:MainProblem:Distance}. We then generate randomly a sequence of 255 feasible candidate solutions $\tilde{x}$, which are used for both \subdirect and \subbenders methods. For each solution, we solve the resulting restricted problem~\eqref{LocalBranching:RestrictedSubproblem} using each method.  Then, we add the Constraints~\eqref{LocalBranching:DistanceOuter} and the Benders optimality cuts derived from the set~\eqref{OptimalityCutParetoMulti} for $\tilde{x}$, and continue to the next solution.

In Table~\ref{TableSubproblemComparison}, we report the average results across all instances in each dataset. The objective value and solve times are as reported as CPLEX, however only instances which terminate before the time limit are included in the solve time. For the subproblems solved it is indicated the percentage of the subproblems which were successfully solved before the time limit was reached.

Examining the results, we note that both methods performed nearly identically in the easy instances. Additionally, in terms of solving time, both methods were very similar in the hard instances. By contrast, the percentage of subproblems solved within the time limit is drastically different between the two methods, with the \subdirect method successfully solving nearly all instances whilst the \subbenders method is able to solve very few. This can also be seen in the objective values, with the average objective value for the \subdirect method being better than the \subbenders method in the hard instances. These results demonstrate quite clearly that the proposed method is much better for solving the subproblem in hard problems compared to the \subbenders method.

A detailed view of the evolution of the solve times in terms of the number of evaluated solutions is presented in Appendix~\ref{AppendixResultsSubproblem}. In brief, we observe a similar increase in the solve time for both methods, with more variability in the \subbenders method.

\begin{table}
\centering
\resizebox{0.9\textwidth}{!}{\begin{tabular}{llrrrrr}
\hline
      &          &    Simple &  Distance &  HomeCharging &   LongSpan &     Price \\
\hline
Objective value & \subdirect &  25382.61 &  12539.70 &      15817.63 &  130526.42 &  32559.12 \\
      & \subbenders &  25382.61 &  12539.70 &      15817.63 &  130462.23 &  32544.02 \\ 
\hline 
Solve time (fully solved, sec) & \subdirect &      0.07 &      0.16 &          0.21 &      30.78 &     14.67 \\
      & \subbenders &      0.08 &      0.14 &          0.20 &      38.84 &     14.39 \\
\hline
Subproblems solved (\%)& \subdirect &     100.00 &      100.00 &          100.00 &      89.04 &     100.00 \\
      & \subbenders &      100.00 &      100.00 &          100.00 &      1.75 &     5.69 \\
\hline
\end{tabular}}
\caption{Comparison of cutting planes solution method and accelerated Benders decomposition for restricted subproblems.}
\label{TableSubproblemComparison}
\end{table}

\subsection{Further results: Local branching subdomain separation scheme}
\label{SectionResultsSpace}

In this section, we compare the two distinct approaches presented in Section~\ref{SectionBranchingScheme} for exploring the disjoint search spaces in the local branching subproblems. We recall that the \branching approach consists in creating a series of branches indicating the time period for which the distance from our candidate exceeds the distance threshold. These branches can be created after any feasible solution, and so we compare branching after every iteration of the Benders decomposition with branching only in iterations of the Benders decomposition for which the solution found is better than the current incumbent. The \disjunctive approach consists in using the Constraints~\eqref{LocalBranching:DistanceOuter} to model the disjunctive set. These constraints are added for every candidate solution encountered during the series of restricted and diversified subproblems.

From  technical standpoint, it is only possible to create two branches at every branch-and-bound node of CPLEX~\citep{CPLEX}. As the \branching approach demands a branch per time period, this comparison is only possible for $T=2$. Since all of the problem instances contain at least four time periods, we modify the instances in this section by not considering any time periods beyond the first two.

In Table~\ref{TableTwoYearSummary}, we see the average results across all instances in each dataset. The solve time, objective value, and optimality gap are reported directly by CPLEX. The number of diversified subproblems and restricted subproblems are collected as part of the callback, and include all subproblems regardless of their solution status. The number of local branching separations is also collected as part of the callback, and reports the number of branches created as part of the separation procedure. 

By examining the results, we observe that the \branching method which branches after every iteration had the worst performance overall. In particular, the objective value is lower and the optimality gap is higher than both other methods. Comparing the \disjunctive and the \branching method which branches only after improving solutions, we see a trade-off in terms of objective value and optimality gap. There is also a notable difference in terms of the number of restricted subproblems, with the \disjunctive method performing roughly 50\% more subproblems. 

To better explain this difference, it is important to keep in mind how both the \branching and \disjunctive methods work at a high level. On the one hand, the \branching method creates two branches at the candidate solution, corresponding to the threshold distance being exceeded in either the first time period or the second. As such, each of these branches solves the main problem~\eqref{LocalBranching:MainProblem} at a local level, within relatively large subdomains. By selecting to explore a subdomain which contains better quality solutions, the solver can then improve the incumbent. However, as the subdomains get more restrictive, the bounds provided by the Benders optimality cuts become more accurate, making it more difficult to find new, improving nodes to explore. Additionally, as a consequence of this branching procedure, any upper bounds that are found are locally valid, leading to a decreased optimality gap in unexplored branches. On the other hand, the \disjunctive method effectively creates a hole within the feasible domain, imposing that all solutions must exceed the threshold distance. As such, the main problem~\eqref{LocalBranching:MainProblem} is working at a global level, not divided into subdomains. This may make it easier to find new, improving nodes, as the solver may easily move to areas of the feasible domain which have not been explored and, thus, areas in which the bounds from the Benders optimality cuts are less accurate.  Additionally, since the main problem works at a global level, any upper bounds that are found are globally valid, leading to an improved optimality gap overall.  Due to the black-box nature of the solver, it is not possible to verify this hypothesis. However, it is supported by the progressively increasing number of nodes and the progressively decreasing optimality gap from the \branching method at all solutions, to the \branching method at improving solutions, and then to the \disjunctive method.

These results suggest that it may be beneficial to consider the \branching method for subproblem separation at improving solutions when higher-quality feasible solutions are sought. However, the limitations in terms of the number of time periods would need to be addressed in order to apply the \branching method more generally. 

Additional performance details are presented in \ref{AppendixResultsSpace}.

\begin{table}
\centering
\resizebox{0.9\textwidth}{!}{\begin{tabular}{llrrrrr}
\hline
                     &                        &    Simple &  Distance &  HomeCharging &  LongSpan &     Price \\
\hline
Solve time (sec) & \branching, all solutions &     46.52 &      2.32 &       1460.34 &   7176.68 &   9045.39 \\
                     & \branching, improving solutions &      0.72 &      1.17 &          8.95 &   7192.53 &   7269.47 \\
                     & \disjunctive &      1.12 &      1.68 &          4.99 &   7175.38 &   7292.31 \\ 
\hline 
Objective value & \branching, all solutions &  13718.14 &   6635.17 &       8151.72 &  18473.34 &  13439.35 \\
                     & \branching, improving solutions &  13718.14 &   6634.57 &       8151.48 &  18882.63 &  13601.62 \\
                     & \disjunctive &  13718.14 &   6634.37 &       8151.48 &  18851.94 &  13552.85 \\ 
\hline 
Optimality gap (\%) & \branching, all solutions &      $<0.01$ &      $<0.01$ &          0.21 &     21.92 &     22.41 \\
                     & \branching, improving solutions &      $<0.01$ &      $<0.01$ &          $<0.01$ &     13.38 &     18.05 \\
                     & \disjunctive &      $<0.01$ &      $<0.01$ &          $<0.01$ &     12.47 &     17.47 \\ 
\hline 
Number of nodes & \branching, all solutions &    509.10 &    100.60 &       3061.90 &   1003.45 &    589.73 \\
                     & \branching, improving solutions &     12.75 &     73.95 &        216.35 &   2654.30 &   1117.95 \\
                     & \disjunctive &     10.00 &     64.10 &        170.80 &   4273.90 &   1592.75 \\ 
\hline 
Number of diversified subproblems & \branching, all solutions &    662.45 &     26.70 &       2892.20 &    760.50 &    278.27 \\
                     & \branching, improving solutions &     13.60 &     13.05 &         61.80 &    889.15 &    276.30 \\
                     & \disjunctive &     18.00 &     15.85 &         30.55 &    876.80 &    268.90 \\ 
\hline 
Number of restricted subproblems & \branching, all solutions &   1424.90 &     98.25 &      13127.65 &   4838.15 &   1708.62 \\
                     & \branching, improving solutions &     31.25 &     41.60 &        195.60 &   4783.15 &   1338.45 \\
                     & \disjunctive  &     87.75 &    101.95 &        170.40 &   6773.35 &   1973.20 \\ 
\hline 
Number of user-created branches & \branching, all solutions &    654.70 &     27.00 &       2887.90 &    741.20 &    260.38 \\
                     & \branching, improving solutions &      2.20 &      4.70 &          6.70 &      12.40 &      9.30 \\
\hline 
\end{tabular}}
\caption{Average performance details for different feasible space reduction methods.}
\label{TableTwoYearSummary}
\end{table}

\section{Conclusion}
\label{SectionConclusion}

In this work, we present two new methods for solving the dynamic MCLP in an exact manner. The accelerated branch-and-Benders-cut method expands upon the current state-of-the-art in the static case, the Benders decomposition in \cite{Cordeau2019}, adding several acceleration techniques from the literature to improve convergence. Computational experiments showed that the proposed method resulted in better performance across solving time, objective value, and optimality gap compared to the state-of-the-art. Additionally, the general nature of the model and the acceleration techniques may allow for the improved method to be applied to different structures of $\Omega$. For example, the dynamic equivalent to the budgeted MCLP~\citep{Khuller1999, Li2021, Wei2023} or the MCLP under uncertainty~\citep{Daskin1983, Berman2013, Vatsa2016, Nelas2020}. 

The accelerated branch-and-Benders-cut method with local branching develops a specialised local branching scheme for the dynamic MCLP. This combines an intuitive distance metric with an innovative subproblem solution method to find improved feasible solutions. Indeed, in computational experiments, this method was the only one able to consistently find better quality feasible solutions compared to the warmstart solution. 

Both methods were applied to an existing problem in the literature, an EV charging station placement model~\citep{Lamontagne2022}. This provided faster solution methods with better performance guarantees, as well as improved lower bounds. These results validate the methodological contributions provided in this work for the dynamic MCLP. 

From our experiments, we conclude that potential speedups should exploit the structure of $\Omega$. In other words, it is crucial to obtain primal formulations in which the linear relaxation is tight in the main problem variables. Hence, future work on specific dynamic MCLPs could follow this research direction in order to improve the performance of our Benders framework.

\section*{Acknowledgements}

The authors gratefully acknowledge the assistance of Jean-Luc Dupr\'e from \emph{Direction Mobilit\'e} of \emph{Hydro-Qu\'ebec} for sharing his expertise on EV charging stations and the network, as well as Ismail Sevim for his insights into the project. We also gratefully acknowledge the assistance of Walter Rei of Universit\'e de Qu\'ebec \`a Montr\'eal and Mathieu Tanneau of Georgia Institute of Technology in discussion about Benders decomposition techniques. 

This research was supported by Hydro-Québec, NSERC Collaborative Research and Development Grant CRDPJ 536757 - 19, and the FRQ-IVADO Research Chair in Data Science for Combinatorial Game Theory.

\clearpage

\input{arxiv.bbl}
\clearpage

\appendix

\section{Additional computational results}
\label{AppendixDetailedResults}

\subsection{Comparison of methodologies}
\label{AppendixResultsBaseline}

In Table~\ref{TableBaselineFull}, we report additional performance details for each of the solution methods compared in Section~\ref{SectionResultsBaseline}. The entries with the solve time, objective value, optimality gap, and the number of nodes have been duplicated in this table for ease of comparison. In addition, we present the number and average time for both lazy cuts (for integer candidate solutions) and user cuts (for fractional candidate solutions). Finally, since the local branching method solves the LP relaxation of the model with the multi-cut method independently, we include the LP solve time here. 

\begin{table}
\centering
\caption{Average performance details for solution methods. }
\label{TableBaselineFull}
\resizebox{0.9\textwidth}{!}{\begin{tabular}{llrrrrr}
\hline
                     &                  &    Simple &  Distance &  HomeCharging &     LongSpan &        Price \\
\hline
Solve time (sec) & Greedy &      0.00 &      0.00 &          0.00 &         0.05 &         0.10 \\
                     & BranchAndCut &      0.20 &      0.70 &          6.56 &      7186.12 &      7194.97 \\
                     & SingleCutBenders &      0.15 &      1.57 &         68.24 &      7184.25 &      7184.76 \\ 
                     & MultiCutBenders &      0.09 &      0.52 &          1.34 &      7184.54 &      7185.04 \\
                     & LocalBranching &      4.67 &      4.81 &         13.00 &      7291.34 &      7305.71 \\

\hline 
Objective value & Greedy &  31814.20 &  16591.75 &      18016.47 &    133724.37 &     33641.74 \\
                     & BranchAndCut &  31820.15 &  16627.14 &      18030.50 &         0.00 &         0.00 \\
                     & SingleCutBenders &  31820.15 &  16627.14 &      18030.47 &    116663.50 &     32563.31 \\                      
                     & MultiCutBenders &  31820.15 &  16627.14 &      18030.40 &    133728.65 &     33642.19 \\
                     & LocalBranching &  31820.15 &  16627.14 &      18030.49 &    133781.88 &     33650.81 \\
\hline 
Optimality gap (\%) & BranchAndCut &      0.00 &      0.00 &          0.00 & - & - \\
                     & SingleCutBenders &      0.00 &      0.00 &          0.00 &  15.79 &        15.00 \\ 
                     & MultiCutBenders &      0.00 &      0.00 &          0.00 &         6.12 &        11.32 \\
                     & LocalBranching &      0.00 &      0.00 &          0.00 &        11.86 &        18.66 \\
\hline 
Number of nodes & BranchAndCut &     59.25 &     43.35 &        131.95 &         0.00 &         0.00 \\
                     & SingleCutBenders &     52.40 &     27.85 &        228.50 &         0.00 &         0.00 \\ 
                     & MultiCutBenders &     65.35 &     75.05 &        860.45 &         0.00 &      1499.60 \\
                     & LocalBranching &     65.05 &    205.60 &        797.35 &      7282.65 &      1202.50 \\
\hline 
Average lazy cut time (sec) & SingleCutBenders &      0.00 &      0.00 &          0.00 &         0.04 &         0.08 \\ 
& MultiCutBenders &      0.00 &      0.00 &          0.00 &         0.06 &         0.12 \\
                     & LocalBranching &      0.20 &      0.31 &          0.51 &       124.45 &        71.14 \\
\hline 
Average user cut time (sec) & SingleCutBenders &      0.00 &      0.00 &          0.00 &         0.05 &         0.11 \\ 
                    & MultiCutBenders &      0.00 &      0.00 &          0.00 &         0.06 &         0.12 \\
                    & LocalBranching &      0.02 &      0.00 &          0.01 &         0.07 &         1.09 \\                  
\hline 
Number of restricted subproblems & LocalBranching &    127.90 &    142.65 &        186.95 &      1042.60 &      1196.00 \\
\hline 
Number of diversified subproblems & LocalBranching &     27.10 &     16.10 &         25.50 &        59.55 &       103.10 \\
\hline 
Solve time, LP (sec) & Multicut &      0.01 &      0.01 &          0.01 &         0.08 &         0.11 \\
\hline
\end{tabular}}
\end{table}

\subsection{Further results: Local branching restricted subproblem solution method}
\label{AppendixResultsSubproblem}

\begin{figure}
    \centering
    \includegraphics[width=0.8\textwidth]{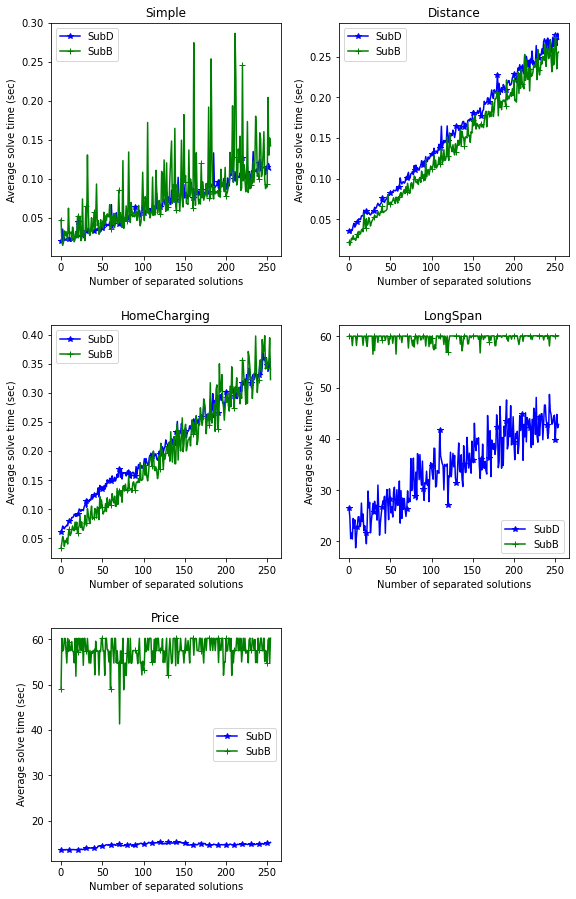}
    \caption{Evolution of subproblem solution time as a factor of the number of subproblem separations completed.}
    \label{FigureSubproblemTime}
\end{figure}

In Figure~\ref{FigureSubproblemTime}, we report the average solve time as a function of the number of separated solutions in each of the datasets. We recall that there are 20 instances in the dataset, and in each instance we randomly generated a series of 255 candidate solutions. The restricted subproblem~\eqref{LocalBranching:RestrictedSubproblem} around each candidate solution was then solved using either the method \textsc{SubD}, which directly solves the reformulation~\eqref{LocalBranching:RestrictedSubproblemReformulated} with a generic solver, or the method \subbenders, which solves the restricted subproblem~\eqref{LocalBranching:RestrictedSubproblem} via the \multicut method. By taking the average solve time of the subproblems across all instances for a fixed number of solved subproblems, we can then examine the impact of the increasing number of cuts to both of the methods.

Looking at Figure~\ref{FigureSubproblemTime}, we observe, for both methods, that the time to solve the restricted subproblem tends to increases along with the number of separated solutions. This is unsurprising, due to the addition of the Big-M Constraints~\eqref{LocalBranching:DistanceOuter:BigM}. More interestingly. we remark a much higher variability in the solving times for the \subbenders method, with the solve times shifting rapidly by 50\% or more. The only dataset for which we see a higher variability in the \subdirect method is in the LongSpan dataset, where nearly all subproblems terminate at the time limit for the \subbenders method.

\subsection{Further results: Local branching subdomain separation scheme}
\label{AppendixResultsSpace}

In Table~\ref{TableTwoYearFull}, we report additional performance details for the subproblem solution methods, as in Section~\ref{SectionResultsSubproblem}. The new values in the table are defined equivalently as for Table~\ref{TableBaselineFull}.

\begin{table}
\centering
\caption{Average performance details for different feasible space reduction methods.}
\label{TableTwoYearFull}
\resizebox{0.9\textwidth}{!}{\begin{tabular}{llrrrrr}
\hline
                     &                        &    Simple &  Distance &  HomeCharging &  LongSpan &     Price \\
\hline
Solve time (sec) & \branching, all solutions &     46.52 &      2.32 &       1460.34 &   7176.68 &   9045.39 \\
                     & \branching, improving solutions &      0.72 &      1.17 &          8.95 &   7192.53 &   7269.47 \\
                     & \disjunctive &      1.12 &      1.68 &          4.99 &   7175.38 &   7292.31 \\ 
\hline 
Objective value & \branching, all solutions &  13718.14 &   6635.17 &       8151.72 &  18473.34 &  13439.35 \\
                     & \branching, improving solutions &  13718.14 &   6634.57 &       8151.48 &  18882.63 &  13601.62 \\
                     & \disjunctive &  13718.14 &   6634.37 &       8151.48 &  18851.94 &  13552.85 \\ 
\hline 
Optimality gap (\%) & \branching, all solutions &      0.00 &      0.00 &          0.21 &     21.92 &     22.41 \\
                     & \branching, improving solutions &      0.00 &      0.00 &          0.00 &     13.38 &     18.05 \\
                     & \disjunctive &      0.00 &      0.00 &          0.00 &     12.47 &     17.47 \\ 
\hline 
Number of nodes & \branching, all solutions &    509.10 &    100.60 &       3061.90 &   1003.45 &    589.73 \\
                     & \branching, improving solutions &     12.75 &     73.95 &        216.35 &   2654.30 &   1117.95 \\
                     & \disjunctive &     10.00 &     64.10 &        170.80 &   4273.90 &   1592.75 \\ 
\hline 
Average lazy cut time (sec) & \branching, all solutions &      0.05 &      0.09 &          0.16 &      8.93 &     33.00 \\
                     & \branching, improving solutions &      0.05 &      0.09 &          0.13 &      9.13 &     31.90 \\
                     & \disjunctive &      0.06 &      0.10 &          0.17 &      9.05 &     33.47 \\ 
\hline 
Number of lazy cuts & \branching, all solutions &      9.50 &      3.25 &         11.65 &     23.85 &     23.15 \\
                     & \branching, improving solutions &      7.90 &      3.35 &         11.00 &     68.90 &     47.50 \\
                     & \disjunctive &     11.45 &     10.05 &          9.40 &     42.10 &     41.35 \\ 
\hline 
Average user cut time (sec) & \branching, all solutions &      0.03 &      0.01 &          0.08 &      3.76 &     12.16 \\
                     & \branching, improving solutions &      0.01 &      0.01 &          0.02 &      1.01 &      3.01 \\
                     & \disjunctive &      0.01 &      0.00 &          0.01 &      0.33 &      0.95 \\ 
\hline 
Number of user cuts & \branching, all solutions &    652.95 &     23.45 &       2880.55 &    736.65 &    255.12 \\
                     & \branching, improving solutions &      5.70 &      9.70 &         50.80 &    820.25 &    228.80 \\
                     & \disjunctive &      6.55 &      5.80 &         21.15 &    834.70 &    227.55 \\ 
\hline 
Number of diversified subproblems & \branching, all solutions &    662.45 &     26.70 &       2892.20 &    760.50 &    278.27 \\
                     & \branching, improving solutions &     13.60 &     13.05 &         61.80 &    889.15 &    276.30 \\
                     & \disjunctive &     18.00 &     15.85 &         30.55 &    876.80 &    268.90 \\ 
\hline 
Number of restricted subproblems & \branching, all solutions &   1424.90 &     98.25 &      13127.65 &   4838.15 &   1708.62 \\
                     & \branching, improving solutions &     31.25 &     41.60 &        195.60 &   4783.15 &   1338.45 \\
                     & \disjunctive &     87.75 &    101.95 &        170.40 &   6773.35 &   1973.20 \\ 
\hline 
Solve time, LP (sec) & \branching, all solutions &      0.00 &      0.01 &          0.01 &      0.05 &      0.24 \\
                     & \branching, improving solutions &      0.00 &      0.00 &          0.00 &      0.05 &      0.27 \\
                     & \disjunctive &      0.00 &      0.01 &          0.01 &      0.14 &      0.62 \\
\hline
Number of user-created branches & \branching, all solutions &    654.70 &     27.00 &       2887.90 &    741.20 &    260.38 \\
                     & \branching, improving solutions &      2.20 &      4.70 &          6.70 &      12.40 &      9.30 \\
\hline 
\end{tabular}}
\end{table}

\end{document}

%% file: Image_SubproblemMethod.tex
\begin{tikzpicture}[Label/.style={text width = 2cm, align = center, font = {\Huge\bfseries\sffamily}},
]

\filldraw[color=blue!80, fill=blue!10, very thick, fill opacity=0.5](0,0) circle (4cm);
\filldraw[color=red!80, fill=red!10, very thick, fill opacity=0.5](-3,-4) circle (4cm);
\filldraw[color=green!80, fill=green!10, very thick, fill opacity=0.5](3,-4) circle (4cm);

\node[Label] (value_j1) at (0,1) {
$10$};
\node[Label] (label_j1) at (3.5,3.5) {
$i_1$};

\node[Label] (value_j2) at (-3.5,-4.5) {
$5$};
\node[Label] (label_j2) at (-6.5,-7.5) {
$i_2$};

\node[Label] (value_j3) at (3.5,-4.5) {
$7$};
\node[Label] (label_j3) at (6.5,-7.5) {
$i_3$};

\node[Label] (intersection_j1j2) at (-1.75,-1.75) {
$8$};
\node[Label] (intersection_j1j3) at (1.75,-1.75) {
$2$};
\node[Label] (intersection_j2j3) at (0,-5) {
$3$};
\node[Label] (intersection) at (0,-3) {
$0$};

\end{tikzpicture}